\documentclass[12pt]{amsart}

\usepackage{amssymb} 
\usepackage{enumerate, amsfonts, latexsym}

\input xy
\xyoption{all}

\newtheorem {theorem}{Theorem} [section]
\newtheorem {lemma} [theorem] {Lemma}

\newtheorem {corollary} [theorem] {Corollary}
\newtheorem {definition} [theorem] {Definition}

\newtheorem{remark} [theorem] {Remark}
\newtheorem{proposition} [theorem] {Proposition}
\newtheorem{convention} [theorem] {Convention}
\newtheorem{terminology} [theorem] {Terminology}
\newtheorem{observation} [theorem] {Observation}

\newfont{\msb}{msbm10 scaled 1200}
\newfont{\euf}{eufm10 scaled 1200}

\def\R {\mathbb R}
\def\F {\mathbb F}
\def\N {\mathbb N}

\def\Z {\mathbb Z}
\def\A {\mathcal A}
\def\B {\mathcal B}

\def\P {\mathbb P}
\def\q {\mathfrak q}

\def\Cinf {\mathcal C_\infty}

\def\SK {\underrightarrow{\text{Ker}}\, }

\begin{document}

\title[Limit groups for relatively hyperbolic groups, I]{Limit groups for relatively hyperbolic groups, I: The basic tools.}

\author[Daniel Groves]{Daniel Groves}
\address{Daniel Groves\\
MSCS UIC 322 SEO, \textsc{M/C} 249\\
851 S. Morgan St.\\
Chicago, IL 60607-7045, USA}
\email{groves@math.uic.edu}

\date{March 20, 2008}


\begin{abstract}
We begin the investigation of $\Gamma$-limit groups, where $\Gamma$ is a torsion-free group which is hyperbolic relative to a collection of free abelian subgroups.  Using the results of \cite{DS}, we adapt the results from \cite{CWIF}.  Specifically, given a finitely generated group $G$, and a sequence of pairwise non-conjugate homomorphisms $\{ h_n : G \to \Gamma \}$, we extract an $\R$-tree with a nontrivial isometric $G$-action.

We then provide an analogue of Sela's shortening argument.
\end{abstract}

\maketitle

In his remarkable series of papers \cite{Sela1, Sela2-6}, Z. Sela has classified those finitely generated groups with the same {\em elementary theory} as the free group of rank $2$ (see also \cite{SelaICM} for a summary).  This class includes all nonabelian free groups, most surface groups, and certain other hyperbolic groups.  In particular, Sela answers in the positive some long-standing questions of Tarski (Kharlampovich and Miasnikov have another approach to these problems; see \cite{KM}).

In \cite{Sela1}, Sela begins with a study of {\em limit groups}.  Sela's definition of a limit group is geometric, though it turns out that a group is a limit group if and only if it is a finitely generated fully-residually free group.  He then produces {\em Makanin-Razborov diagrams}, which give a parametrization of $\text{Hom}(G,\mathbb F)$, where $G$ is an arbitrary finitely generated group and $\mathbb F$ is a nonabelian free group (such a parametrisation is also given in \cite{KM2}).  Over the course of his six papers, two of the main tools Sela uses are the theory of isometric actions on $\R$-trees and the {\em shortening argument}.  

Sela's work naturally raises the question of which other classes of groups can be understood using this approach.  Many of Sela's methods (and, more strikingly, some of the answers) come from geometric group theory.  Thus, when looking for classes of groups to apply these methods to, it seems natural to consider groups of interest in geometric group theory.  In \cite{SelaHyp}, Sela considers an arbitrary torsion-free hyperbolic group $\Gamma$, and characterises those groups with the same elementary theory as $\Gamma$.  Of particular note is his results that any group which has the same elementary theory as a torsion-free hyperbolic group is itself a torsion-free hyperbolic group.  This result exhibits a deep connection between the logic of groups and geometric group theory.  In \cite{Alibegovic2}, Alibegovi\'c constructs Makanin-Razborov diagrams for limit groups.  In \cite{CWIF} the author began this study for certain torsion-free CAT$(0)$ groups.

This paper serves two purposes.  First, we generalise the results of \cite{CWIF} to the context of torsion-free groups which are hyperbolic relative to a collection of abelian subgroups\footnote{See Section \ref{RelHyp} below for a definition and discussion of relatively hyperbolic groups.}.  We
construct a space closely related to the Cayler graph (see Section \ref{ConstructX}), and use the results of Dru\c{t}u and Sapir
from \cite{DS} to analyse an asymptotic cone of this space.  We then follow \cite{CWIF} to extract an $\R$-tree from this asymptotic cone.\footnote{The $\R$-tree we construct is different to that constructed by Dru\c{t}u and Sapir in \cite{DS-trees}.  See Remark \ref{WhyThisPaper?} below for a discussion of the reasons we feel the $\R$-tree in this paper is still a worthwhile construction in light of the seemingly more general construction in \cite{DS-trees}.}  Armed with this $\R$-tree, we then develop an analogue of Sela's shortening argument in this context.

The following result is a straightforward application of the shortening argument (see Section \ref{ShortenSection} for a definition of $\text{Mod}(\Gamma)$).

\begin{theorem} \label{ModinAut}
Suppose that $\Gamma$ is a torsion-free group which is hyperbolic relative to a collection of free abelian subgoups.  Then ${\rm Mod}(\Gamma)$ has finite index in ${\rm Aut}(\Gamma)$.
\end{theorem}

The true context of this paper is as the beginning of the study of $\Gamma$-limit groups where $\Gamma$ is torsion-free and hyperbolic relative to free abelian subgroups.  In a continuation paper \cite{MR-RH} we use the results of this paper and of \cite{SelaHyp} to construct {\em Makanin-Razborov diagrams} for such a group $\Gamma$.  It is our hope that much, possibly all, of Sela's program can be carried out for these groups.

The outline of this paper is as follows: In Section \ref{LimitSection} we recall the concepts of {\em limit groups} and {\em $\Gamma$-limit groups}.  In Section \ref{RelHyp} we recall the definition of relatively hyperbolic groups, and the basic results required for this paper.  In Section \ref{DSResults} we recall the concept of asymptotic cones and some results of Dru\c{t}u and Sapir from \cite{DS}.  In Section \ref{ConstructX} we define a space $X$, closely related to the Cayley graph of the relatively hyperbolic group $\Gamma$.  The space $X$, equipped with a natural $\Gamma$-action, seems to be an appropriate space  with which to study $\Gamma$-limit groups.  In Section \ref{Asymp} a particular asymptotic cone $X_\omega$ is extracted from a sequence of non-conjugate homomorphisms $\{ h_n : G \to \Gamma \}$, where $G$ is an arbitrary finitely generated group and the limiting action of $G$ on $X_\omega$ is studied.  In Section \ref{RTree} we extract an action of $G$ on an $\R$-tree with no global fixed point.  In Section \ref{ShortenSection} we state a version of Sela's shortening argument in the context of this paper. In Section \ref{ActionsonRtrees} we recall some of the theory of groups acting isometrically on $\R$-trees.
In Sections \ref{IET}--\ref{Discrete} we present the proof of Theorem \ref{ShorteningArgument}
(The Shortening Argument).

\begin{remark} \label{WhyThisPaper?}
This paper subsumes the results which were in the first version of this paper, and also in
\cite{CWIF2}.  The exception to this is the proof that the groups under consideration are
Hopfian.  This result follows immediately from \cite[Theorem 5.2]{MR-RH} (and the proof
there does not depend on anything left out of this version of this paper).  Since the proof of
the Hopf property is technical, and no more enlightening than the proof of \cite[Theorem 5.2]{MR-RH}, we chose to leave this result out of this paper.  Leaving this proof out made it natural to take the results
of \cite{CWIF2} and the first version of this paper and merge them into this single paper.

Subsequent to this paper appearing on the arxiv, the paper \cite{DS-trees} of Dru\c{t}u and Sapir
appeared.  In this, they take a finitely generated group $G$ acting on an arbitrary tree-graded metric space and produce an isometric action of $G$ on an $\R$-tree.  

In this paper, we require that the tree-graded metric space has pieces isometric to Euclidean
spaces, and that the stabilizers of these pieces act by translations.  However, there are benefits
to the construction in this paper over that of \cite{DS-trees}.  Namely, as we see in this paper and in \cite{MR-RH}, with our tree we are able to prove the appropriate analogue of Sela's shortening argument.  Thus we feel that our construction is still important.  

The construction in this paper has been used in \cite{MR-RH} and \cite{Dahmani-G}.  It will also
be crucial in future work of the author on the elementary theory of relatively hyperbolic groups.
\end{remark}

{\bf Acknowledgements.}
{\em The paper \cite{CWIF} was written before \cite{DS} appeared.  Although I had already had a cursory glance at \cite{DS}, I would like to thank both Mark Sapir and Jon McCammond for suggesting that the results in \cite{DS} may allow the results from \cite{CWIF} to be generalised as they have been in this paper.  I would also like to thank Jason Manning and Indira Chatterji for fruitful discussions on the construction in Section \ref{ConstructX}.

The first version of this paper was written whilst the author was a Taussky-Todd Instructor at Caltech.  The author was supported in part by the NSF.  Both of these organisations are thanked for their support.}

\section{Limit groups and $\Gamma$-limit groups} \label{LimitSection}

Recall the following two definitions, due to Bestvina and Feighn \cite{BFSela}.

\begin{definition} \cite[Definition 1.5]{BFSela}
Let $G$ and $\Xi$ be finitely generated groups, and let $\{ h_n : G \to \Xi \}$ be a sequence of homomorphisms.  The {\em stable kernel} of $\{ h_n \}$, denoted $\SK(h_n)$, is the set of all elements $g \in G$ so that $g \in \text{ker}(h_n)$ for all but finitely many $n$.

The sequence $\{ h_n \}$ is {\em stable} if for all $g \in G$, either (i) $g \in \SK(h_n)$; or (ii) $g \not\in \text{ker}(h_n)$ for all but finitely many $n$.  
\end{definition}

\begin{definition} \cite[Definition 1.5]{BFSela} \label{DefAlgLimit}
A {\em $\Xi$-limit group} is a group of the form $G/\SK(h_n)$ where $G$ is a finitely generated group and $\{ h_n : G \to \Xi \}$ is a stable sequence of homomorphisms.
\end{definition}

\begin{remark}
If each of the $h_n$ is equal to $h$, a single homomorphism, then the sequence $\{ h_n \}$ is certainly stable and the associated $\Xi$-limit group is just $h(G)$.  In particular, all finitely generated subgroups of $\Xi$ are $\Xi$-limit groups.
\end{remark}

A {\em limit group} is an $\F$-limit group, where $\F$ is a finitely generated free group.  This terminology is due to Sela \cite{Sela1}, although the definition that Sela gave was in terms on an action of $G$ on an $\R$-tree induced by the sequence $\{ h_n : G \to \F \}$.  Sela's geometric definition also makes sense for $\delta$-hyperbolic groups (see \cite{SelaHyp}).  In this paper, we pursue a geometric definition of $\Gamma$-limit groups, where $\Gamma$ is a torsion-free group which is hyperbolic relative to free abelian subgroups.

In case $\Xi = \F$, the geometric and algebraic definitions of $\Xi$-limit groups are the same (see \cite[Lemma 1.3]{Sela1}).   The two definitions are also the same when $\Xi$ is a torsion-free $\delta$-hyperbolic group (see \cite[Lemma 1.3]{SelaHyp}).  When $\Xi$ is a torsion-free CAT$(0)$
 group with isolated flats whose flat stabilisers are abelian, a geometric definition of $\Xi$-limit group was given in \cite[Definition 3.21]{CWIF} and it was proved \cite[Theorem 5.1]{CWIF} that these two definitions are the same.

Suppose that $\Gamma$ is a torsion-free group hyperbolic relative to free abelian subgroups.  In this paper, we provide an appropriate geometric definition of $\Gamma$-limit group, in analogy with the definition from \cite{CWIF} (see Definition \ref{DefGeoLimit} below).  It is proved in Theorem \ref{TwoDefsSame} that this definition is equivalent to Definition \ref{DefAlgLimit}.  As in Sela's definition, along with the geometric definition comes a faithful action of a (strict) $\Gamma$-limit group on an $\R$-tree.

The utility of the algebraic Definition \ref{DefAlgLimit} is that it has implications for the logic of $\Gamma$.  In the case of Sela's limit groups, the nonabelian limit groups are exactly those that have the same universal theory as a nonabelian free group.  In general, if $T_\forall(H)$ denotes the {\em universal theory} of a group $H$ then we have the following (see \cite{CG} for a detailed discussion of this issue)

\begin{lemma}
Let $\Xi$ be a finitely presented group and suppose that $L$ is a $\Xi$-limit group.  Then $T_\forall(\Xi) \subseteq T_\forall(L)$.
\end{lemma}

The utility of Sela's geometric definition is that it allows the application of the (Rips) theory of isometric actions on $\R$-trees, and Sela uses this to make a very deep study of limit groups (and of $\Gamma$-limit groups, where $\Gamma$ is a torsion-free hyperbolic group).  It turns out that the class of limit groups is exactly the class of fully residually free groups, which has been widely studied in the past.

\section{Relatively hyperbolic groups} \label{RelHyp}

Recently there has been a large amount of interest in relatively hyperbolic groups.  Relatively hyperbolic groups were originally defined by Gromov in his seminal paper \cite{Gromov}, and an alternative definition was given by Farb \cite{Farb}.  Bowditch \cite{Bowditch} gave two definitions, equivalent to Gromov's and Farb's, respectively (see \cite{Dahmani} for a proof of the equivalence of the definitions).  Dru\c{t}u and Sapir \cite{DS} gave a characterisation of relatively hyperbolic groups in terms of their asymptotic cones.  The results of this paper rely heavily on the results of \cite{DS}.

Examples of relatively hyperbolic groups include: (i) geometrically finite Kleinian groups (which are hyperbolic relative to their cusp subgroups); (ii) fundamental groups of hyperbolic manifolds of finite volume (hyperbolic relative to their cusp subgroups); (iii) hyperbolic groups (relative to the trivial group, or a finite collection of quasi-convex subgroups); (iv) free products (relative to the factors); and (v) limit groups (relative to their maximal noncyclic abelian subgroups).  See \cite{Farb, Bowditch, Dahmani2, Sz} for details. 

For further recent work on relatively hyperbolic groups, see \cite{Alibegovic, Yaman, DS, DS--RD, DS-trees, Dahmani3, G-Manning, Osin, CR} (among others).

The definition of relatively hyperbolic which we give is a hybrid of Farb's definition and a definition of Bowditch.

\begin{definition}[Coned-off Cayley graph]
Suppose that $\Gamma$ is a finitely generated group, with finite generating set $\A$, and that $\{ H_1, \ldots , H_m \}$ is a collection of finitely generated subgroups of $\Gamma$.  Let $X$ be the Cayley graph of $\Gamma$ with respect to $\A$.  We form the {\em coned-off Calyey graph}, $\tilde{X}$, by adding to $X$ a vertex $c_{\gamma,H_i}$ for each coset $\gamma H_i$ of a parabolic subgroup, and for each coset $\gamma H_i$, an edge from $c_{\gamma,H_i}$ to $\gamma'$ for each $\gamma' \in \gamma H_i$.
\end{definition}

\begin{definition}
We say that $\Gamma$ is {\em hyperbolic relative to $\{ H_1, \ldots , H_m \}$} if 
\begin{enumerate}
\item the coned-off Cayley graph $\tilde{X}$ is $\delta$-hyperbolic for some $\delta$; and
\item for each edge $e \in \tilde{X}$, and each $n \ge 1$, there are only finitely many loops of length at most $n$ which contain $e$.
\end{enumerate}
\end{definition}

\begin{terminology}
Suppose that $\Gamma$ is a group which is hyperbolic relative to the collection $\{ H_1, \ldots , H_m \}$ of subgroups.  The subgroups $H_i$ are called {\em parabolic subgroups}.\footnote{Alternative terminology for these subgroups is {\em peripheral} subgroups.  Sometimes, all conjugates of the $H_i$ are also called parabolic subgroups.}
\end{terminology}

In this paper we are concerned with torsion-free relatively hyperbolic groups $\Gamma$ where all the parabolic subgroups are free abelian.

\begin{definition}
A subgroup $K$ of a group $G$ is {\em malnormal} if for all $g \in G \smallsetminus K$ we have $gKg^{-1} \cap K = \{ 1 \}$.

A group $G$ is {\em CSA} if any maximal abelian subgroup of $G$ is malnormal.
\end{definition}

\begin{lemma}  
Suppose that $\Gamma$ is a torsion-free group which is hyperbolic relative to a collection of free abelian subgroups.  Then $\Gamma$ is CSA.
\end{lemma}
\begin{proof}
Let $A$ be a maximal abelian subgroup of $\Gamma$.

Since $\Gamma$ is torsion-free, the Bounded Coset Penetration property implies that any conjugate of a parabolic subgroup is malnormal (see \cite[Example 1, p.819]{Farb}).  This implies that if $M$ is a conjugate of a parabolic subgroup and $A$ intersects $M$ nontrivially then $A = M$, which is malnormal.

Suppose that $A$ is a maximal abelian subgroup of $\Gamma$ and that $g \in A \smallsetminus \{ 1 \}$. If $g$ is not contained in a conjugate of a parabolic subgroup then a result of Osin (see \cite[Theorem 1.14, p.10]{Osin} and the comment thereafter) implies that the centraliser of $\langle g \rangle$ is virtually cyclic.  Since $\Gamma$ is torsion-free, this centraliser is cyclic.  Therefore, in this case $A = \langle h \rangle$ for some $h$.  Note that $A$ is maximal cyclic in $\Gamma$.  Suppose now that $\gamma \in \Gamma$ is such that $\gamma h^k \gamma^{-1} = h^j$ for some $k,j \in \Z \smallsetminus \{ 0 \}$.  Then \cite[Corollary 4.21, p.83]{Osin} implies that $|k| = |j|$.  Thus, $\gamma^2$ commutes with $h^j$. This implies that $\gamma^2 \in \langle h \rangle$, which in turn implies that $\gamma \in \langle h \rangle$, so $A$ is malnormal.
\end{proof}

\section{Dru\c{t}u and Sapir's results} \label{DSResults}

In \cite{DS}, Dru\c{t}u and Sapir find a characterisation of relatively hyperbolic groups in terms of their asymptotic cones.  In this section, we recall the definition of asymptotic cones and then briefly summarise those of Dru\c{t}u and Sapir's results necessary for this paper.

\subsection{Asymptotic cones}

Asymptotic cones were introduced by van den Dries and Wilkie in \cite{VW} in order to recast and simplify Gromov's Polynomial Growth Theorem from \cite{GromovPoly}.  See \cite{DS} for a discussion of other results about asymptotic cones.  We briefly recall the definition of asymptotic cones.

\begin{definition}
A {\em non-principal ultrafilter}, $\omega$, is a $\{ 0,1 \}$-valued finitely additive measure on $\mathbb N$ defined on all subsets of $\mathbb N$ so that any finite set has measure $0$.
\end{definition}
The existence of non-principal ultrafilters is guaranteed by Zorn's Lemma.  We fix once and for all a non-principal ultrafilter $\omega$.\footnote{A different choice of ultrafilter can change the resulting asymptotic cone in interesting ways, but not in a way that affects our results.  Thus, we are unconcerned which ultrafilter is chosen.}  Given any bounded sequence $\{ a_n \} \subset \R$ there is a unique number $a \in \R$ so that for all $\epsilon > 0$ we have $\omega( \{ a_n \mid |a-a_n| < \epsilon \} ) = 1$.  We denote $a$ by $\omega$-$\lim \{ a_n \}$.  This notion of limit exhibits most of the properties of the usual limit (see \cite{VW}).

Let $(X,d)$ be a metric space.  Suppose that $\{ \mu_n \}$ is a sequence of real numbers with no bounded subsequence, and that $\{ x_n \}$ is a collection of points in $X$.  Let $(X_n,d_n)$ be the metric space which has set $X$ and metric $\frac{1}{\mu_n}d_X$.  The {\em asymptotic cone of $X$ with respect to $\{ x_n \}, \{ \mu_n \}$ and $\omega$}, denoted $X_\omega$, is defined as follows.  First, define the set $\tilde{X_\omega}$ to consist of all sequences $\{ y_n \mid y_n \in X_n \}$ for which $\{ d_{X_n}(x_n, y_n) \}$ is a bounded sequence.  Define a pseudo-metric $\tilde{d}$ on $\tilde{X_\omega}$ by
\[	\tilde{d}(\{ y_n \}, \{ z_n \}) = \mbox{$\omega$-$\lim$} \{ d_{X_n}(y_n,z_n) \}	.	\]
The asymptotic cone $X_\omega$ is the metric space induced by the pseudo-metric $\tilde{d}$ on $\tilde{X_\omega}$:
\[	 X_\omega := \tilde{X_\omega} / \sim		,	\]
where the equivalence relation $`\sim '$ on $\tilde{X_\omega}$ is defined by: $x \sim y$ if and only if $\tilde{d}(x,y) = 0$.  The pseudo-metric $\tilde{d}$ on $\tilde{X_\omega}$ naturally descends to a metric on $d_\omega$ on $X_\omega$.

\subsection{Tree-graded spaces and relatively hyperbolic spaces}

\begin{definition} \cite[Definition 1.10]{DS} \label{TreeGradedDef}
Let $Y$ be a complete geodesic metric space, and let $\mathcal P$ be a collection of closed subsets of $Y$ (called {\em pieces}).  We say that the space $Y$ is {\em tree-graded with respect to $\mathcal P$} if the following two conditions are satisfied:
\begin{enumerate}
\item[$(T_1)$] Each pair of distinct pieces intersect in at most a point; and
\item[$(T_2)$] Every simple geodesic triangle (a simple loop composed of three geodesics) in $X$ is contained in a single piece.
\end{enumerate}
\end{definition}

It is worth remarking that in \cite{CWIF} it is proved that if $Y$ is a CAT$(0)$ space with isolated flats and relatively thin triangles then a particular asymptotic cone $Y_\omega$ of $Y$ is tree-graded with respect to its collection of maximal flats (the proof of this is essentially contained in \cite{KL}).  It was this fact that inspired the current paper.  In \cite{HK}, Hruska and Kleiner prove that if a cocompact CAT$(0)$ space has isolated flats then it has relatively thin triangles.

One of the main results of \cite{DS} is the following

\begin{theorem} \cite[Theorem 1.11]{DS}
A finitely generated group $G$ is relatively hyperbolic with respect to finitely generated subgroups $H_1, \ldots , H_n$ if and only if every asymptotic cone of $G$ (with respect to any non-principal ultrafilter, any sequence of scaling constants, where the basepoints are the identity of $G$) is tree-graded with respect to $\omega$-limits of sequences of cosets of the subgroups $H_i$.
\end{theorem}

We also need the following definition and results.

\begin{definition} \cite[Definition 3.19]{DS}
Let $(Y,\text{dist})$ be a metric space and let $\mathcal Q = \{ Q_i \mid i \in I \}$ be a collection of subsets of $Y$.  In every asymptotic cone $Y_\omega$, with choice of basepoints $\{ x_n \}$, we consider the collection of subsets
\[	\mathcal Q_\omega = \left\{ \text{lim}^\omega (Q_{i_n}) \mid (i_n)^\omega \in I^\omega \mbox{ such that $\left\{ \frac{\text{dist}(x_n,Q_{i_n})}{d_n} \right\}$ is bounded} \right\} .	\]
The space $Y$ is {\em asymptotically tree-graded with respect to $\mathcal Q$} if every asymptotic cone $Y_\omega$ is tree-graded with respect to $\mathcal Q_\omega$.
\end{definition}

\begin{theorem} \cite[Theorem 5.1]{DS}
Let $Y$ be a metric space and let $\mathcal Q$ be a collection of subsets of $Y$.  Let $\mathfrak q : Y \to Y'$ be a quasi-isometry.  If $Y$ is asymptotically tree-graded with respect to $\mathcal Q$ then $Y'$ is asymptotically tree-graded with respect to $\mathfrak q (\mathcal Q)$.
\end{theorem}

\begin{theorem} \cite[Theorem 4.1]{DS} \label{DS41}
Let $(Y,\text{dist})$ be a geodesic metric space and let $\mathcal Q = \{ Q_i \mid i \in I \}$ be a collection of subsets of $Y$.  The space $Y$ is asymptotically tree-graded with respect to $\mathcal Q$ if and only if the following properties are satisfied:
\begin{enumerate}
\item[$(\alpha_1)$] For every $\xi > 0$, the diameters of the intersections $\mathcal N_\xi(Q_i) \cap \mathcal N_\xi(Q_j)$ are uniformly bounded for all $i \neq j$;
\item[$(\alpha_2)$] For every $\theta \in [0,\frac{1}{2})$ there exists $M(\theta) > 0$ so that for every geodesic $\q$ of length $l$ and every $Q \in \mathcal Q$ with $\q (0), \q (l) \in \mathcal N_{\theta l}(Q)$ we have $\q ([0,l]) \cap \mathcal N_M(Q) \neq \emptyset$;
\item[$(\alpha_3)$] For every $k \ge 2$ there exists $\zeta > 0$, $\nu \ge 8$ and $\chi > 0$ such that every $k$-gon $P$ in $X$ with geodesic edges which is $(\zeta, \nu, \chi)$-fat satisfies $P \subseteq \mathcal N_\chi(Q)$ for some $Q \in \mathcal Q$.
\end{enumerate}
\end{theorem}

\section{The space $X$} \label{ConstructX}

In this paragraph we find a space $X$, closely associated to the Cayley graph of a relatively hyperbolic group, which will be the appropriate space for our analysis of $\Gamma$-limit groups in the subsequent sections, and also in \cite{MR-RH} and \cite{Dahmani-G}.

Suppose that $\Gamma$ is a group which is hyperbolic relative to a collection $\{ H_1, \ldots , H_m \}$ of subgroups.  

Choose a generating set $\A$ for $\Gamma$ which intersects each of the subgroups $H_i$ in a generating set $\mathcal B_i$ for $H_i$, for $1 \le i \le m$. Let $\mathcal B = \cup_{i=1}^m \mathcal B_i$.  For each $i \in \{ 1 ,\ldots , m \}$.

Let $d_{\A}$ be the word metric on $\Gamma$ induced by the generating set $\A$ of $\Gamma$, and let $d_{\mathcal B_i}$ be the word metric on $H_i$ induced by $\mathcal B_i$.

Let $Y$ denote the Cayley graph of $\Gamma$ with respect to $\A$, where each edge is isometric to the unit interval $[0,1]$.  The group $\Gamma$ acts on itself by left multiplication, which induces an isometric action on $Y$.

Let $\gamma H_i$ be a coset of some parabolic subgroup of $\Gamma$.  The set $\mathcal B_i$ also gives a metric on $\gamma H_i$, which we denote by $d_{\mathcal B_i}$.  Now, \cite[Lemma 4.3]{DS} states that there is a constant $K \ge 0$ so that for any $x, y \in \gamma H_i$, any geodesic joining $x$ and $y$ in $Y$ stays entirely in the $K$-neighbourhood of $\gamma H_i$.

We now build a space $Y^k$ out of $Y$.  Indira Chatterji told me of a similar construction which she attirbuted to David Epstein.  Epstein proved an analogue of Theorem \ref{CuspIsomEmbed} for his space.

Consider a coset $\gamma H_i$, along with the set of edges labelled by elements of $\mathcal B_i$.  The resulting subgraph $Z(\gamma ,H_i)$ of $Y$ is exactly the Cayley graph of $H_i$.  We now form a new graph $Z(\gamma H_i)^1$, which is another copy of $Z(\gamma ,H_i)$, except that each edge is isometric to the closed interval $[0,\frac{1}{4}]$.  Denote this new graph by $Z(\gamma,H_i)^1$, and join it to $Z(\gamma,H_i)$ by joining corresponding vertices by edges of length $\frac{1}{4}$.  Perform this construction for each coset $\gamma H_i$ of a parabolic subgroup.

We define $Y^j$ inductively, starting from $Y^{j-1}$.  First, form $Z(\gamma,H_i)^j$ with edges of length $2^{-2j}$ (but otherwise isomorphic to $Z(\gamma,H_i)$), and join it to $Z(\gamma,H_i)^{j-1}$ by edges of length $2^{-2j}$.   

\begin{terminology}
We call the edges of length $2^{-2j}$ joining
$Z(\gamma,H_i)^{j-1}$ to $Z(\gamma,H_i)^j$ {\em vertical}, and the edges which lie in
some $Z(\gamma,H_i)^j$ {\em horizontal}.
\end{terminology}
The space $Y^j$ is the union of $Y^{j-1}$ along with $Z(\gamma,H_i)^j$, and the vertical edges joining $Z(\gamma,H_i)^{j-1}$ to $Z(\gamma,H_i)^j$.  Endow $Y^j$ with the natural path metric.

For an integer $k \ge 1$, let $C(\gamma,H_i)^k$ be the union of the graphs $Z(\gamma,H_i), Z(\gamma,H_i)^1, \ldots , Z(\gamma,H_i)^k$, along with the sets of vertical edges that join successive graphs in this sequence.

There is a natural space $Y^{\infty}$, the metric completion of $\cup_{s=1}^\infty Y^s$, (where we consider $Y^s$ to be a subset of $Y^{s+1}$).  Each coset $\gamma H_i$ inherits a `cone-point' $w_{\gamma,i}$ from this completion process.  In $Y^\infty$, the point $w_{\gamma,i}$ lies at distance $\eta$ from the coset $\gamma H_i$, where $\eta = \sum_{s=1}^\infty 2^{-2s} < \frac{1}{2}$.  It is clear that the space $Y^{\infty}$ is quasi-isometric to the coned-off Cayley graph of $\Gamma$.  Let $Y^{\infty}$ be $\Upsilon$-hyperbolic, and suppose without loss of generality that $\Upsilon > 1$.

\begin{definition}
Suppose that $H_i$ is a parabolic subgroup of $\Gamma$, and $\gamma H_i$ is a coset of $H_i$ in $\Gamma$.  The space $P_{\gamma ,i}$ is formed from $Y$ by performing the construction of $Y^\infty$ to all cosets of parabolic subgroups {\em except} the coset $\gamma H_i$.  
\end{definition}

\begin{lemma} \label{qconvex}
Suppose that $\gamma H_i$ is a coset of a parabolic subgroup in $\Gamma$, and that $x,y \in \gamma H_i$.  Let $[x,y]$ be a geodesic between $x$ and $y$ in $P_{\gamma,i}$ and suppose that $[x,y]$ does not intersect $\gamma H_i$ except at its endpoints.  Then $[x,y]$ lies entirely in the $35 \Upsilon$-neighbourhood of $\gamma H_i$ in $P_{\gamma,i}$.
\end{lemma}
\begin{proof}
Let $\widehat{[x,y]}$ be the image of $[x,y]$ in the space $Y^\infty$, under the inclusion $P_{\gamma,i} \subset Y^\infty$.  For any $R$, the $R$-neighbourhood of $w_{\gamma,i}$ in $Y^\infty$ naturally corresponds to the $(R- \eta)$-neighbourhood of $\gamma H_i$ in $P_{\gamma,i}$, and if $R > \eta$ then outside of these balls the two spaces are locally isometric (where the local isometry is induced by the inclusion $P_{\gamma,i} \subset Y^\infty$).

Suppose that $\widehat{[x,y]}$ is not contained in the $10\Upsilon$ neighbourhood of $w_{\gamma,i}$ in $Y^\infty$.  That part of $\widehat{[x,y]}$ which lies at least $10\Upsilon$ from $w_{\gamma,i}$ is a $10\Upsilon$-local-geodesic.  By \cite[Theorem III.H.1.13, p. 405]{BH}, outside the $10\Upsilon$ ball around $w_{\gamma,i}$ in $Y^\infty$, the path $\widehat{[x,y]}$ is a $(\frac{14\Upsilon}{6\Upsilon},2\Upsilon)$-quasi-geodesic.  However, the distance in $Y^\infty$ which it travels outside of the $10\Upsilon$ ball around $w_{\gamma,i}$ is at most $20\Upsilon$ (since the path starts and finishes at distance $\eta < \frac{1}{2}$ from $w_{\gamma,i}$).

Therefore, the total distance that $\widehat{[x,y]}$ travels outside the $10\Upsilon$ ball about $w_{\gamma ,i}$ is at most
\[	\left(\frac{7}{3}\right)20\Upsilon + 2 \Upsilon < 50\Upsilon .	\]
Therefore $\widehat{[x,y]}$ is contained in the $35\Upsilon$ ball around $w_{\gamma,i}$ in $Y^\infty$.  As above, this implies that $[x,y]$ is contained in the $35\Upsilon$-neighbourhood of $\gamma H_i$ in $P_{\gamma,i}$, as required.
\end{proof}

\begin{lemma} \label{qisomembed}
There exists a constant $K_1$, depending only on $Y$, and the set $\{ \gamma H_i \}$, so that for all $x, y \in \gamma H_i$
\[	d_{P_{\gamma,i}}(x,y) \le d_{\gamma H_i}(x,y) \le K_1 d_{P_{\gamma,i}}(x,y).	\]
\end{lemma}
\begin{proof}
Since $\gamma H_i \subset P_{\gamma,i}$, and since $d_{P_{\gamma,i}}$ is a path metric, the first inequality is immediate.

Let $x, y \in \gamma H_i$, and let $[x,y]$ be a geodesic between $x$ and $y$ in $P_{\gamma,i}$.  By Lemma \ref{qconvex} above, the path $[x,y]$ lies entirely within the $35\Upsilon$-neighbourhood of $\gamma H_i$.  

Let $c_1, c_2, \ldots , c_k$ be points along $[x,y]$ which are such that $\eta \le d_{P_{\gamma,i}}(c_i,c_{i+1}) \le 1$ (this can be ensured if we choose each $c_i$ to be either a vertex from $Y$ or a cone-point $w_{\gamma',j}$).  For each $1 \le i \le k$, let $b_i \in \gamma H_i$ be a point in $\gamma H_i$ as close as possible to $c_i$, and choose a geodesic $[c_i,b_i]$ (which has length at most $35\Upsilon$).  Possibly $c_i = b_i$, and so the path $[c_i,b_i]$ is a constant path  Also, choose a path $[b_i,b_{i+1}] \subset \gamma H_i$ of shortest length.

Consider the paths $p_i = [b_i,b_{i+1}]$ and $q_i = [b_i,c_i,c_{i+1},b_{i+1}]$.  The path $q_i$ has length at most $70\Upsilon +1$.  Also, unless $[c_i,c_{i+1}] \subset \gamma H_i$, the path $q_i$ intersects $\gamma H_i$ only at its endpoints.  

The path $q_i$ corresponds to a path $q_i' \subset Y$, where any part of $q_i$ which passes through a cone-point is replaced by a (shortest) path through the corresponding coset.  Now, $q_i'$ is a relative $(70\Upsilon+1)$-quasi-geodesic.  Also, $p_i$ is a relative $2\eta$-quasi-geodesic, and can be considered as a path in $Y$.  Note that $p_i$ penetrates $\gamma H_i$ while $q_i'$ does not.  Therefore, Bounded Coset Penetration implies that there is a constant $c = c(70 \Upsilon +1)$  so that $p_i$ travels distance at most $c$ in $\gamma H_i$, which is to say that $p_i$ has length at most $c$.

We have seen that each $[b_i,b_{i+1}]$ has length at most $c$.  Therefore $d_{\gamma H_i}(x,y) \le c k$.  However, $d_{P_{\gamma,i}}(x,y) \ge \eta k$, and it suffices to take $K_1 = \frac{c}{\eta}$.
\end{proof}

Now let $P_{\gamma,i}^{(k)}$ be the space formed from $Y$ by adding the spaces $C(\gamma',H_j)^k$ for all cosets of parabolic subgroups except $\gamma H_i$.  Then for all $x,y \in \gamma H_i$ we have
\[	d_{P_{\gamma,i}}(x,y) \le d_{P_{\gamma,i}^{(k)}}(x,y) 	,	\]
and so by Lemma \ref{qisomembed} we have
\[	d_{\gamma H_i}(x,y) \le K_1 d_{P_{\gamma,i}^{(k)}}(x,y)	.	\]

\begin{theorem} \label{CuspIsomEmbed}
There exists $k \ge 0$ so that each of the graphs $Z(\gamma,H_i)^k$ is isometrically embedded in $Y^k$.
\end{theorem}
\begin{proof}
It suffices to take $k > \frac{\log_2 K_1}{2}$.

Suppose that there exists $u,v \in Z(\gamma,H_i)^k$ so that a geodesic $[u,v]$ between $u$ and $v$ does not lie entirely within $Z(\gamma,H_i)^k$.  Since $Z(\gamma,H_i)^k$ is isometrically embedded in $C(\gamma,H_i)^k$, the path $[u,v]$ cannot be contained in $C(\gamma,H_i)^k$.  

Suppose that $x$ is the point furthest along $[u,v]$ so that $[u,x] \subset C(\gamma,H_i)^k$.  Let $y$ be the point furthest along $[u,v]$ so that $[x,y]$ intersects $C(\gamma,H_i)^k$ only in $Z(\gamma,H_i)$.  Because of the way $C(\gamma,H_i)^k$ was built, that part of $[u,v]$ immediately before $x$ consists entirely of edges joining the different $Z(\gamma,H_i)^i$, from $Z(\gamma,H_i)^k$ to $Z(\gamma,H_i)$.  Similarly, that part of $[u,v]$ immediately after $y$ consists of a `vertical' path from $y$ to $Z(\gamma,H_i)$.  Let $x_1$ be the final point in $[u,x]$ contained in $Z(\gamma,H_i)^k$ and let $y_1$ be the first point in $[y,v]$ contained in $Z(\gamma,H_i)^k$.  Then we have $d_{\gamma H_i}(x,y) = 2^{2k}d_{Z(\gamma,H_i)^k}(x_1,y_1)$.  Now let $D = d_{P_{\gamma,i}^{(k)}}(x,y)$, and note that $D \ge \eta$.
\begin{eqnarray*}
D + \frac{1}{2} & \le & d_{Z(\gamma,H_i)^k}(x_1,y_1) \\
& = & 2^{-2k} d_{\gamma H_i}(x,y) \\
& \le & 2^{-2k} K_1 D.
\end{eqnarray*}
Therefore, $D+\frac{1}{2} \le 2^{-2k}K_1 D$, which implies in particular that $2^{-2k}K_1 -1 > 0$, contradicting our choice of $k$.  This completes the proof.
\end{proof}

We now fix $k$ so that Theorem \ref{CuspIsomEmbed} holds, and consider the space $Y^k$.

\begin{lemma} \label{CosetsAttract}
There exists a function $f_1 : \N \to \N$ so that if $x, y \in Y^k$ are such $x$ and $y$ lie in the $N$-neighbourhood of $\gamma H_i$ and $[x,y]$ does not intersect $Z(\gamma,H_i)^k$ then $d(x,y) \le f_1(N)$.
\end{lemma}
\begin{proof}
By Theorem \ref{CuspIsomEmbed}, and the definition of $Y^k$, it suffices to bound the length of a geodesic $[w,z]$ where $w, z \in \mathcal N_N(\gamma H_i)$ and $[w,z]$ does not intersect $C(\gamma,H_i)^k \smallsetminus Z(\gamma,H_i)$.

For such a pair we have $d_{Y^k}(w,z) = d_{P_{\gamma,i}^{(k)}}(w,z)$.  Denote this distance by $E$.  Let $w_1, z_1$ be points in $\gamma H_i$ which are closest to $w$ and $z$, respectively.  Also, let $w_2,z_2$ be the points in $Z(\gamma,H_i)^k$ which are closest to $w_1$ and $z_1$, respectively.  Also, let $\eta_k = \sum_{i=1}^k 2^{-2i}$ be the distance from $\gamma H_i$ to $Z(\gamma,H_i)^k$.  Then we have
\begin{eqnarray*}
E & = & d_{Y^k}(w,z)\\
& \le & d_{Y^k}(w_1,z_1) + 2N\\
& = & d_{Z(\gamma,H_i)^k}(w_2,z_2) + 2N + 2\eta_k\\
& = & 2^{-2k} d_{\gamma H_i}(w_1,z_1) + 2N + 2\eta_k\\
& = & 2^{-2k} K_1 d_{P_{\gamma,i}^{(k)}}(w_1,z_1) + 2N + 2\eta_k\\
& \le & 2^{-2k} K_1 \left( d_{P_{\gamma,i}^{(k)}}(w,z) + 2N \right) + 2N + 2\eta_k\\
& = & 2^{-2k}K_1 \left( E + 2N \right) + 2N + 2\eta_k ,
\end{eqnarray*}
which implies (since the choice of $k$ from Theorem \ref{CuspIsomEmbed} ensures that $1 - 2^{-2k}K_1 > 0$) that
\[	E \le \frac{2^{-2k}K_1 N + 2N + 2\eta_k}{1 - 2^{-2k} K_1} 	.	\]
This completes the proof.
\end{proof}

We now assume that $\Gamma$ is torsion-free and that each of the parabolic subgroups of $\Gamma$ are free abelian.  

\begin{remark}
Suppose that a group $G$ is hyperbolic relative to a family $\mathcal P$ of subgroups, and that
some of the subgroups in $\mathcal P$ are hyperbolic.  Let $\mathcal P'$ be the non-hyperbolic
groups in $\mathcal P$.  Then $G$ is also hyperbolic relative to $\mathcal P'$.

Therefore, we assume that all parabolic subgroups of our relatively hyperbolic groups are not
hyperbolic.  In case parabolics are free abelian, as in this paper, this amounts to assuming that
none of the parabolics are infinite cyclic (or trivial).
\end{remark}

Suppose that the generating set for $\Gamma$ intersects each
parabolic subgroup in a basis (as a free abelian group).  Then each of the graphs $Z(\gamma,H_i)^k$ is isomorphic to the `standard' Cayley graph of $\Z^{n_i}$ (with edge of length $2^{-2k}$.  Fix an embedding $\phi_i : Z(\gamma,H_i)^k \hookrightarrow \R^{n_i}$ which is isometric on
each edge, and send the vertices adjacent to the identity to (scaled) standard basis vectors
in $\R^{i_k}$ (and their negatives).

Using the map $\phi_i$, glue a copy of $\R^{n_i}$ onto each subspace $Z(\gamma,H_i)^k$ of $Y^k$ where $n_i$ is the rank of $H_i$, and $\R^{n_i}$ is equipped with the standard Euclidean ($\ell_2$-) metric.

\begin{definition}
The resulting space is denoted $X$, and $\mathcal Q$ is the collection of copies of  the $\R^{n_i}$ glued onto the cosets $\gamma H_i$ (where $i \in \{ 1, \ldots , m \}$ and $\gamma \in \Gamma$). 
\end{definition}

The copies of $\R^n$ that have been glued to $Y^k$ to form the space $X$ now play the role of cosets.  They are isometrically embedded, and Lemma \ref{CosetsAttract} above holds for these subspaces also, since lengths of paths are unchanged outside of $Z(\gamma,H_i)^k$, and distances can only get shorter inside $Z(\gamma,H_i)^k$.

The action of $\Gamma$ on $X$ is defined in the obvious way.  The stabiliser in $\Gamma$ of any $Q \in \mathcal Q$ is a conjugate of a parabolic subgroup, which acts by translations on $Q$.

\subsection{Properties of $X$}

\begin{lemma}
Suppose $Q \in \mathcal Q$ is a copy of $\R^{n_i}$ in $X$ as above.  Then $Q$ is isometrically embedded and convex in $X$.
\end{lemma}

Given this lemma, we call the elements $Q \in \mathcal Q$ `flats'.

\begin{lemma}
Left multiplication of $\Gamma$ on itself induces an isometric action of $\Gamma$ on $X$.  This action is proper and cocompact.
\end{lemma}

\begin{lemma} \label{XAsymTree}
$X$ is asymptotically tree-graded with respect to the set $\mathcal Q$.
\end{lemma}
\begin{proof}
The inclusion map $\Gamma \hookrightarrow X$ is a quasi-isometry.  

We know that $\Gamma$ is asymptotically tree-graded with respect to the set of cosets $\gamma H_i$.  
By the proof of \cite[Theorem 5.1, p.44]{DS}, $X$ is asymptotically tree-graded with respect to $\mathcal Q$.
\end{proof}

Note also that any asymptotic cone of $\Gamma$ is bi-Lipschitz homeomorphic to the analogous asymptotic cone of $X$ (taking the same basepoints, and the same scaling factors).  The utility of using $X$ rather than just $\Gamma$ is the following

\begin{lemma}
Suppose that $X_\omega$ is an asymptotic cone of $X$.  Then each piece of $X_\omega$ is isometric to $(\R^k, d)$ for some $k$, where $d$ is the standard Euclidean metric on $\R^k$.
\end{lemma}

In analogy with Hruska's definition of {\em Isolated Flats} for CAT$(0)$ spaces (see \cite[2.1.2]{Hruska}), we note the following

\begin{lemma} [Isolated Flats] \label{IsolatedFlats}
Let $\mathcal Q$ be the collection of flats in $X$.  Then there is a function $\phi : \R_+ \to \R_+$ such that for every pair of distinct flats $Q_1, Q_2 \in \mathcal Q$ and for every $k \ge 0$, the intersection of the $k$-neighbourhoods of $Q_1$ and $Q_2$ has diameter less than $\phi(k)$.
\end{lemma}
\begin{proof}
This is merely a restatement of Theorem \ref{DS41}.($\alpha_1$).
\end{proof}

\begin{convention} \label{IsolConv}
Let $\phi : \R_+ \to \R_+$ be as in Lemma \ref{IsolatedFlats}.  We suppose that $\phi(k) \ge k$ for all $k \ge 0$ and that $\phi$ is a nondecreasing function.
\end{convention}

We now prove a quasi-convexity result for the metric on $X$

\begin{lemma} \label{MetricQConvex}
There exists a function $N_1 ; \N \to \N$ so that for any $K$, if $x_1,x_2,y \in X$ so that $d_X(x_1,x_2) \le  K$ and $[x_1,y]$ and $[x_2,y]$ are geodesics, then $[x_1,y]$ is contained in the $N_1(K)$-neighbourhood of $[x_2,y]$ (and vice versa).
\end{lemma}
\begin{proof}
Choose a geodesic $[x_1,x_2]$.  Then the path $[y,x_1,x_2] = [y,x_1] \cup [x_1,x_2]$ is a $(1, K)$-quasi-geodesic.

By \cite[Theorem 1.12]{DS} there are constants $\tau$ and $M$ so that:

$\bullet$ $[y,x_1,x_2]$ is contained in the $\tau$-tubular neighbourhood of the $M$-saturation of $[y,x_2]$ (see \cite[Definition 8.9]{DS} for a definition of $M$-saturations); and

$\bullet$  the points at which $[y,x_1,x_2]$ enters and leaves the $\tau$-neighbourhood of flats in the $M$-saturation of $[y,x_2]$ are at bounded distance from $[a,x_1]$.

By Lemma \ref{CosetsAttract}, and the fact the flats are isometric to $\R^n$, the path $[y,x_1,x_2]$ lies in the $D_2$-neighbourhood of $[y,x_2]$ for some constant $D_2$.  A symmetric argument on $[y,x_2,x_1]$ and $[y,x_1]$ implies that $[y,x_2,x_1]$ lies in the $D_2$-neighbourhood of $[y,x_1]$.
\end{proof}

For our purposes, one of the most important properties of the space $X$ is contained in the following theorem, which shows that geodesic triangles in $X$ satisfy the {\em Relatively Thin Triangles Property} (see \cite[Definition 3.1.1]{Hruska}).

\begin{theorem} \label{RelThinTriangles}
Suppose that $X$ is as constructed above.  There exists $\delta > 0$ so that for any $a,b,c \in X$, and any $\Delta(a,b,c)$ is a geodesic triangle, either (i) $\Delta(a,b,c)$ is $\delta$-thin in the usual sense; or else (ii) there is a unique flat $E \subset X$ so that each side of $\Delta(a,b,c)$ is contained in the $\delta$-neighbourhood of the union of $E$ and the other two sides.
\end{theorem}
\begin{proof}
Choose a geodesic triangle $\Delta(a,b,c)$ in $X$, with a choice of geodesics $[a,b]$, $[b,c]$ and $[a,c]$.

By \cite[Lemma 8.16]{DS} and \cite[Lemma 8.17]{DS}, there is a constant $\alpha$ (independent of the points $a,b,c$) such that one of two possibilities occurs: either 

(i) there is a point $x \in X$ whose $\alpha$-neighbourhood intersects all three of the geodesics $[a,b]$, $[b,c]$ and $[a,c]$ nontrivially; or

(ii) there is a flat $E \in \mathcal Q$ so that the $\alpha$-neighbourhood of $E$ intersects each of the three geodesics nontrivially.

In case (i), let $x_1$ be a point on $[a,b]$ which is within $\alpha$ of $x$, and let $x_2$ be a point on $[a,c]$ which is within $\alpha$ of $x$.  Then $d_X(x_1, x_2) \le 2\alpha$.

In case (ii) let $x_1$ be the point on $[a,b]$ which is closest to $a$ subject to being in the $\alpha$-neighbourhood of $E$, and similarly for $x_2$ on $[a,c]$.  Then \cite[Corollary 8.14]{DS} implies that there is a constant $D_1$ so that $d_X(x_1,x_2) \le D_1$.  We assume that $D_1 \ge 2\alpha$.

Therefore, in either case, there exist points $x_1 \in [a,b]$ and $x_2 \in [a,c]$ so that $d_X(x_1,x_2) \le D_1$.  Denote by $[a,x_1]$ the sub-path of $[a,b]$ from $a$ to $x_1$, and similarly for $[a,x_2] \subset [a,c]$.  By Lemma \ref{MetricQConvex}, there is a constant $D_2$ so that $[a,x_1]$ lies in the $D_2$-neighbourhood of $[a,x_2]$, and vice versa.

We use a symmetric argument on the points $b$ and $c$ -- finding points $y_1 \in [c,a]$, $y_2 \in [c,b]$ and $z_1 \in [b,a]$, $z_2 \in [b,c]$ as with $x_1$ and $x_2$.

Now, in case (i) above, we can take $x_1 = z_1$, $x_2 = y_1$ and $y_2 = z_2$, and we're done.  In case (ii), we note that the path $[x_1, z_1] \subseteq [a,b]$ lies in the $N_1(\alpha)$neighbourhood of $E$, by Lemma \ref{MetricQConvex}, and similarly for $[x_2,y_1] \subseteq [a,c]$ and $[z_2,y_2] \subseteq [b,c]$.  

Therefore, it suffices to take $\delta = \max \{ D_2, N_1(\alpha) \}$.
\end{proof}

\subsection{Projecting to flats}

In this paragraph we record some results about projecting to flats which are required for the proofs in the subsequent sections.

\begin{definition} \cite[Definition 4.9]{DS}
Let $x \in X$ and $A \subset X$.  The {\em almost projection of $x$ onto $A$} is the set of points $y \in A$ so that $d_X(x,y) \le d_X(x,A) +1$.
\end{definition}

The following result follows immediately from \cite[Corollary 8.14]{DS} and Theorem \ref{XAsymTree}.

\begin{lemma}
There exists a constant $C_1$ so that if $Q \in \mathcal Q$ and $x \in X$ then the almost projection of $x$ onto $Q$ has diameter at most $C_1$.
\end{lemma}

The following result also follows immediately from \cite[Corollary 8.14]{DS} and Theorem \ref{XAsymTree}.

\begin{lemma} \label{ProjectConvex}
There exists a function $N_3 : \N \to \N$ so that if $x_1, x_2 \in X$, $Q \in \mathcal Q$ and $\pi(x_1), \pi(x_2)$ are in the almost projections of $x_1$ and $x_2$ to $Q$, respectively then $d_X(\pi(x_1),\pi(x_2)) \le N_3(d_X(x_1,x_2))$.
\end{lemma}

Again, we suppose that $N_3(x) \ge x$ for all $x \ge 0$ and that $N_3$ is a nondecreasing function.

Recall that $\delta$ is the constant from Theorem \ref{RelThinTriangles} and that $\phi : \N \to \N$ is the function from Lemma \ref{IsolatedFlats}.

\begin{lemma} [cf. Lemma 2.11, \cite{CWIF}] \label{UniqueFlat}
Suppose that $\Delta = \Delta(a,b,c)$ is a geodesic triangle in $X$.  If $\Delta$ is not $\left( \delta + \frac{\phi(\delta)}{2} \right)$-thin then $\Delta$ is $\delta$-thin relative to a {\em unique} flat $Q \in \mathcal Q$.
\end{lemma}
\begin{proof} Given Lemma \ref{IsolatedFlats} and Theorem \ref{RelThinTriangles}, the proof of \cite[Lemma 2.11]{CWIF} applies directly.
\end{proof}

\begin{lemma} [cf. Lemma 2.21, \cite{CWIF}] \label{222}
Suppose that $Q \in \mathcal Q$, that $x,y \in Q$ and that $z \in X$.  Let $[x,z]$ and $[y,z]$ be geodesics.  Then there exist $u \in [x,z]$ and $v \in [y,z]$ that both lie in the $2\delta$-neighbourhood of $Q$ such that
\[	d_X(u,v) \le \phi(\delta)	.	\]
\end{lemma}
\begin{proof}
Given Theorems \ref{CuspIsomEmbed} and \ref{RelThinTriangles} and Lemma \ref{IsolatedFlats}, the proof of \cite[Lemma 2.21]{CWIF} applies directly.
\end{proof}

\begin{proposition} [cf. Proposition 2.22, \cite{CWIF}] \label{ProjectProp}
Suppose that $Q \in \mathcal Q$, that $x,y \in X$ and that some geodesic $[x,y]$ does not intersect the $4\delta$-neighbourhood of $Q$.  Let $\pi(x) , \pi(y)$ be in the almost projections of $x$ and $y$, respectively.  Then 
\[	d_X(\pi(x), \pi(y)) \le   N_3(\phi(3\delta) + N_1(\phi(\delta))).	\]
\end{proposition}
\begin{proof}
By Lemma \ref{222} there exist $w_1 \in [\pi(x),y]$ and $w_2 \in [\pi(y),y]$, both in the $2\delta$-neighbourhood of $Q$ such that $d_X(w_1,w_2) \le \phi(\delta)$.  By a similar argument as in the proof of Lemma \ref{222} (see \cite{CWIF}), there are $u_1 \in [\pi(x),x]$ and $u_2 \in [\pi(x),y]$ which lie outside the $2\delta$-neighbourhood of $E$ so that $d_X(u_1,u_2) \le \phi(3\delta)$. Now $[w_1,y]$ is contained in the $N_1(\phi(\delta))$-neighbourhood of $[w_2,y]$, by Lemma \ref{MetricQConvex}.  Therefore, there exists $u_3 \in [\pi(y),y]$ so that $d_X(u_2,u_3) \le N_1(\phi(3\delta))$.

We can choose $\pi(u_1)$ and $\pi(u_3)$ in the almost projections of $u_1,u_3$ so that $\pi(u_1) = \pi(x)$ and $\pi(u_3) = \pi(y)$. Now, $d_X(u_1,u_3) \le \phi(3\delta) + N_1(\phi(\delta))$, so by Lemma \ref{ProjectConvex} we have
\begin{eqnarray*}
d_X(\pi(x),\pi(y)) & = & d_X(\pi(u_1),\pi(u_3)) \\
& \le & N_3(d_X(u_1,u_3)) \\
& \le & N_3(\phi(3\delta) + N_1(\phi(\delta))),	
\end{eqnarray*}
as required.
\end{proof}

\section{Asymptotic cones and compactification} \label{Asymp}

In this section we start with $\Gamma$, a finitely generated group which acts properly and cocompactly by isometries on a metric space $(X,d_X)$, a finitely generated group $G$ and a sequence $\{ h_n : \Gamma \to G \}$ of homomorphisms.  Using $\{ h_n \}$ we construct a particular asymptotic cone $X_\omega$, which is equipped with an isometric action of $G$ with no global fixed point.

In the case of $\delta$-hyperbolic groups and spaces, the construction we describe in this section is essentially due to Paulin \cite{Paulin, Paulin2} (see also Bestvina \cite{Bestvina} and Bridson--Swarup \cite{BS}), though was not cast there in terms of asymptotic cones.  For CAT$(0)$ spaces, this construction is performed by Kapovich and Leeb \cite{KL}.  The general construction is similar. See \cite{CWIF} for more details about this construction and \cite{VW} or \cite{DS} for many properties about asymptotic cones.

Let $G$ be a finitely generated group and $\Gamma$ a torsion-free group which is hyperbolic relative to a collection of free abelian subgroups.  Let $\A$ be a finite generating set for $G$, let $X$ be the space constructed from a Cayley graph of $\Gamma$ in Section \ref{ConstructX}, and let $x \in X$ correspond to the identity of $\Gamma$.  If $h : G \to \Gamma$ is a homomorphism, define
\[	\| h \| := \min_{\gamma \in \Gamma}\max_{g \in \A} d_X(x, (\gamma h(g) \gamma^{-1}) . x)	,	\]
and let $\gamma_h$ be an element of $\Gamma$ which realises this minimum.

\begin{terminology}
We say that a pair of homomorphisms $h, h' : G \to \Gamma$ are {\em non-conjugate} if there is no inner automorphism $\tau : \Gamma \to \Gamma$ so that $h' = \tau \circ h$.
\end{terminology}

Suppose that $\{ h_i : G \to \Gamma \}$ is a sequence of pairwise non-conjugate homomorphisms.  Then the sequence $\{ \| h_n \| \}$ does not contain a bounded subsequence.  Let $X_\omega$  be the asymptotic cone, defined with respect to some non-principal ultrafilter $\omega$, the sequence of basepoints $x_n = x$ and the sequence of scaling factors $\mu_n = \| h_n \|$.

The action of $G$ on $X_\omega$ is defined by $g . \{ y_n \} = \{ \gamma_{h_n} h_n(g) \gamma_{h_n}^{-1} . y_n \}$.

\begin{lemma} \label{NoFixedPt}
The action of $G$ on $X_\omega$ has no global fixed point.
\end{lemma}
\begin{proof}
See \cite[Lemma 3.9]{CWIF}, where the proof does not use the CAT$(0)$ property.
\end{proof}

\subsection{The action of $G$ on $X_\omega$} \label{GAction}

\begin{definition}
Define a separable $G$-invariant subspace $\Cinf$ of $X_\omega$ to be the union of (i) the geodesic segments $[x_\omega, g . x_\omega]$ for all $g \in G$; and (ii) the flats $Q_\omega \subseteq X_\omega$ which contain a simple geodesic triangle contained in $\Delta (g_1 . x_\omega, g_2 . x_\omega , g_3 . x_\omega )$ for some $g_1 , g_2 , g_3$.
\end{definition}

\begin{lemma}
The space $\Cinf$ is (i) separable; (ii) $G$-invariant; (iii) convex in $X_\omega$; and (iv) tree-graded with pieces isometric to $(\R^n, \ell_2)$, for some $n$ (which may vary according to the piece).
\end{lemma}

Suppose that $\{ (Y_n,\lambda_n) \}_{n=1}^\infty$ and $(Y,\lambda)$ are pairs consisting of metric spaces, together with actions $\lambda_n : G \to \text{Isom}(Y_n)$, $\lambda : G \to \text{Isom}(Y)$.  Recall (cf. \cite[\S 3.4, p. 16]{BFSela}) that $(Y_n, \lambda_n ) \to (Y,\lambda)$ in the {\em $G$-equivariant Gromov topology} if and only if:  for any finite subset $K$ of $Y$, any $\epsilon > 0$ and any finite subset $P$ of $G$, for sufficiently large $n$, there are subsets $K_n$ of $Y_n$ and bijections $\rho_n : K_n \to K$ such that for all $s_n, t_n \in K_n$ and all $g_1, g_2 \in P$ we have
\[	\left| d_{Y}(\lambda(g_1) . \rho_n (s_n) , \lambda(g_2) . \rho_n (t_n)) - d_{Y_n} ( \lambda_n(g_1) . s_n ,  \lambda_n(g_2) . t_n ) \right| < \epsilon .	\]

To a homomorphism $h : G \to \Gamma$, we naturally associate a pair $(X_h,\lambda_h)$ as follows: let $X_h = X$, endowed with the metric $\frac{1}{\mu_h}d_X$; and let $\lambda_h = \iota \circ h$, where $\iota: \Gamma \to \text{Isom}(X)$ is the fixed homomorphism.

Let $\lambda_\infty : G \to \text{Isom}(\Cinf)$ denote the action of $G$ on $\Cinf$.

\begin{proposition} \cite[Lemma 3.15]{CWIF} \label{GromovTop}
If there is a separable $G$-invariant subspace $\mathcal C$ of $X_\omega$ which contains the basepoint $x_\omega$ of $X_\omega$  then there is a subsequence $\{ f_i \}$ of $\{ h_i \}$ so that $(X_{f_i},\lambda_{f_i}) \to (\Cinf, \lambda_\infty)$ in the $G$-equivariant Gromov topology.
\end{proposition}

For the remainder of the section and the next, we assume that we have passed to the convergent subsequence $\{ f_i \}$ of $\{ h_i \}$.  In this vein, we denote $X_{f_i}$ by $X_i$, and $\lambda_{f_i}$ by $\lambda_i$.

\begin{lemma} \cite[Corollary 3.17]{CWIF} \label{Flats}
Let $\mathcal F_\infty$ be the set of flats in $\mathcal C_\infty$.  For each $E \in \mathcal F_\infty$ there is a sequence $\{ E_i \subset X_i \}$ so that $E_i \to E$ in the $G$-equivariant Gromov topology.
\end{lemma}

\begin{observation}
The action of $G$ on $\mathcal C_\infty$ has no global fixed point.
\end{observation}

\begin{definition} \label{DefGeoLimit}
Suppose that $G$ and $\Gamma$ are finitely generated groups and that $\{ h_i : G \to \Gamma \}$ is a sequence of pairwise non-conjugate homomorphisms, leading to an isometric action of $G$ on $\Cinf$, where $\Cinf$ is constructed from $X_\omega$, the asymptotic cone of $\Gamma$, as above.  Let $K_\infty$ be the kernel of the action of $G$ on $\Cinf$:
\[	K_\infty = \{ g \in G \ | \ \forall y \in \Cinf, g . y = y \} .	\]
The {\em strict $\Gamma$-limit group} is $L_\infty = G/ K_\infty$.

A {\em $\Gamma$-limit group} is a group which is either a strict $\Gamma$-limit group as above or else a finitely generated subgroup of $\Gamma$.
\end{definition}

The following result is clear from the definition of the Gromov topology.

\begin{lemma} \label{SKinKinf}
Suppose that the sequence of homomorphisms $\{ f_i : G \to \Gamma \}$ gives rise to a sequence of actions converging to an action of $G$ on $\Cinf$, and that $K_\infty$ is the kernel of the action of $G$ on $\Cinf$.  Then $\SK (f_i) \subset K_\infty$.
\end{lemma}

The following results give information about the flats in $\mathcal C_\infty$, and their stabilisers in $G$.

\begin{proposition} [cf. \cite{CWIF}, Lemma 3.18] \label{FlatInv}
Suppose $g \in G$ leaves a flat $E \subseteq \Cinf$ (setwise) invariant, and that $\{ E_j \}$ converges to $E$.  Then for all but finitely many $i$ we have $f_i(g) . E_i = E_i$.
\end{proposition}
\begin{proof} 
The proof of \cite[Lemma 3.18]{CWIF} applies directly.
\end{proof}

\begin{proposition} [cf. \cite{CWIF}, Lemma 3.19] \label{ActTrans}
Suppose $g \in \text{Stab}_G(E)$ for some flat $E \subseteq \Cinf$.  Then $g$ acts (possibly trivially) by translation on $E$.
\end{proposition}
\begin{proof}
The proof of \cite[Lemma 3.19]{CWIF} applies directly, once we notice that an element of $\gamma \in \Gamma$ which leaves a flat in $X$ invariant lies in a conjugate of a parabolic subgroup and acts by Euclidean translations on the flat.
\end{proof}

\section{The $\R$-tree $T$} \label{RTree}

\subsection{Constructing the $\R$-tree} \label{ConstructT}

We now follow \cite{CWIF} to construct from $\Cinf$ an $\R$-tree $T$ equipped with an isometric $G$-action with no global fixed point.  Given the construction of $\mathcal C_\infty$ in the previous section, the construction of $T$ is exactly the same as in \cite{CWIF}.  We repeat the definition of $T$ here.

Let $\mathcal F_\infty$ be the collection of all pieces in $\mathcal C_\infty$.  By Definition \ref{TreeGradedDef}, for any $g \in G$ exactly one of the following holds: (i) $g . E = E$; (ii) $|g . E \cap E | = 1$; or (iii) $g . E \cap E = \emptyset$.  By Lemmas \ref{FlatInv} and \ref{ActTrans}, $\text{Stab}(E)$ is a countable abelian group, acting by translations on $E$ (possibly not faithfully).  

Let $\mathcal D_E$ be the set of directions of the translations of $E$ by elements of $\text{Stab}(E)$.

For each element $g \in G \smallsetminus \text{Stab}(E)$, let $l_g(E)$ be the (unique) point where any geodesic from a point in $E$ to a point in $g . E$ leaves $E$, and let $\mathcal L_E$ be the set of all $l_g(E) \subset E$.  Note that if $g.E \cap E$ is nonempty (and $g \not\in \text{Stab}(E)$) then $g . E \cap E = \{ l_g(E) \}$.

Since $G$ is finitely generated, and hence countable, both sets $\mathcal D_E$ and $\mathcal L_E$ are countable.  Given a (straight) line $p \subset E$, let $\chi_E^p$ be the projection from $E$ to $p$.  Since $\mathcal L_E$ is countable, there is a line $p_E \subset E$ such that:
\begin{enumerate}
\item the direction of $p_E$ is not orthogonal to a direction in $\mathcal D_E$; and
\item if $x$ and $y$ are distinct points in $\mathcal L_E$, then $\chi_E^{p_E}(x) \neq \chi_E^{p_E}(y)$.
\end{enumerate}
Project $E$ onto $p_E$ using $\chi_E^{p_E}$.  The action of $\text{Stab}(E)$ on $p_E$ is defined in the obvious way (using projection) -- this is an action since the action of $\text{Stab}(E)$ on $E$ is by translations.  Connect $\mathcal C_\infty \smallsetminus E$ to $p_E$ in the obvious way -- this uses the following

\begin{observation}
Suppose $S$ is a component of $\mathcal C_\infty \smallsetminus E$.  Then there is a (unique) point $x_S \in E$ so that $S$ is a component of $\mathcal C_\infty \smallsetminus \{ x_S \}$.
\end{observation}

Glue such a component $S$ to $p_E$ at the point $\chi_E^{p_E}(x_S)$.

Perform this projecting and gluing construction in an equivariant way for all flats $E \subseteq \mathcal C_\infty$ -- so that for all $E \subseteq \mathcal C_\infty$ and all $g \in G$ the direction of the lines $p_{g . E}$ and $g . p_E$ is the same (this is possible since the action of $\text{Stab}(E)$ on $E$ is by translations, so doesn't change directions).

Having done this for all flats $E \subseteq \mathcal C_\infty$, we arrive at a space $T$ which we endow with the (obvious) path metric.

An isometric action of $G$ on $T$ is naturally induced by the action of $G$ on $X_\omega$.

The space $T$ has a distinguished set of geodesic lines, namely those of the form $p_E$, for $E \in \mathcal F_\infty$.  Denote the set of such geodesic lines by $\mathbb P$.

The following lemma is \cite[Lemma 4.2]{CWIF}, and the proof there holds in the current situation.

\begin{lemma} \label{NoFixedPtonT}
The space $T$ is an $\R$-tree and has an isometric $G$-action with no global fixed point.
\end{lemma}

\begin{remark}
Since $K_\infty \le G$ acts trivially on $\mathcal C_\infty$, it also acts trivially on $T$, and the action of $G$ on $T$ induces an isometric action of $L_\infty$ on $T$.
\end{remark}

\subsection{The actions of $G$ and $L_\infty$ on $T$}

Let $G$ be a finitely generated group, and $\Gamma$ a torsion-free group which is hyperbolic relative to a collection of free abelian subgroups.  Suppose that $\{ h_i : G \to \Gamma \}$ is a sequence of pairwise non-conjugate homomorphisms.  Let $X_\omega$, $\Cinf$ and $T$ be as in Section \ref{Asymp} and Subsections \ref{GAction} and \ref{ConstructT}, repsectively.  Let $\{ f_i : G \to \Gamma \}$ be the subsequence of $\{ h_i \}$ as in the conclusion of Proposition \ref{GromovTop}.  Let $K_\infty$ be the kernel of the action of $G$ on $\Cinf$ and let $L_\infty = G/ K_\infty$ be the associated strict $\Gamma$-limit group.

\begin{theorem} \label{LinfProps} [cf.  \cite{Sela1}, Lemma 1.3; \cite{CWIF}, Theorem 4.4]
In the above situation, the following properties hold.
\begin{enumerate}
\item Suppose that $[A,B]$ is a nondegenerate segment in $T$.  Then $\text{Stab}_{L_\infty}[A,B]$ is an abelian subgroup of $L_\infty$; \label{SegStabAb}
\item If $T$ is isometric to a real line then for all but finitely many $n$ the group $f_n(G)$ is free abelian.  Furthermore, in this case $L_\infty$ is free abelian; \label{Tline}
\item If $g \in G$ fixes a tripod in $T$ pointwise then $g \in \SK (f_i )$; \label{tripod}
\item Let $[y_1, y_2] \subset [y_3, y_4]$ be a pair of non-degenerate segments of $T$ and assume that $\text{Stab}_{L_\infty}[y_3,y_4]$ is non-trivial.  Then
\[	\text{Stab}_{L_\infty}[y_1,y_2] = \text{Stab}_{L_\infty}[y_3,y_4] .	\]
In particular, the action of $L_\infty$ on the $\R$-tree $T$ is stable; \label{Stable}
\item Let $g \in G \smallsetminus K_\infty$.  Then for all but finitely many $n$ we have $g \not\in \text{ker}(f_n)$; \label{seqStable}
\item $L_\infty$ is torsion-free; \label{tf} and
\item If $T$ is not isometric to a real line then $\{ f_i \}$ is a stable sequence of homomorphisms. \label{fiStable}
\end{enumerate}
\end{theorem}
\begin{proof}
The proof of \cite[Theorem 4.4]{CWIF} relies on a number of different results.  In each case, we have an exact analogue in the setting of Theorem \ref{LinfProps} here.

The results we need are:  Proposition \ref{FlatInv}, Lemma \ref{SKinKinf}, Lemma \ref{Flats}, Lemma \ref{UniqueFlat}, Lemma \ref{222}, Proposition \ref{GromovTop}, Proposition \ref{ProjectProp} and the fact that stabilisers in $\Gamma$ of flats in $X$ are malnormal (see\cite[Example 1, p. 819]{Farb}).

Given these results, the proof of \cite[Theorem 4.4]{CWIF} applies directly.\footnote{Note in particular
that the analogue of \cite[Lemma 4.5]{CWIF} holds in this setting.  We state this separately below.}
The only change is that some of the constants have changed, so some of the counting has to be changed.  This is straightforward.
\end{proof}

We need the following lemma later when we describe the shortening argument.

\begin{lemma} \cite[Lemma 4.5]{CWIF} \label{4.5}
Suppose $\alpha, \beta \in X$ and $g \in G$ are such that there is a segment of length at 
least $$6 \phi(4\delta) + 4 \max\{ d_X(g\alpha,\alpha), d_X(g\beta, \beta) \}$$ in a geodesic
$[\alpha, \beta]$ which is within $\delta$ of a flat $E \in \mathcal{Q}$.  Then $g \in \rm{Fix}(E)$.
\end{lemma}
\begin{proof}
Given Theorem \ref{RelThinTriangles} and Lemma \ref{IsolatedFlats}, the proof of \cite[Lemma 4.5]{CWIF} applies without change.
\end{proof}

The following are two immediate applications of the above construction of the $\R$-tree $T$, and of Theorem \ref{LinfProps}.  See \cite{CWIF} for proofs which apply without change in the current setting.

\begin{theorem} [cf. Theorem 5.1, \cite{CWIF}] \label{TwoDefsSame}
Suppose that $\Gamma$ is a torsion-free group which is hyperbolic relative to a collection of free abelian subgroups.  A group $L$ is a $\Gamma$-limit group in the sense of Definition \ref{DefAlgLimit} if and only if it is a $\Gamma$-limit group in the sense of Definition \ref{DefGeoLimit}.
\end{theorem}

\begin{theorem} [cf. Theorem 5.9, \cite{CWIF}] \label{SplittingTheorem}
Suppose that $\Gamma$ is a torsion-free group which is hyperbolic relative to a collection of free abelian subgroups, and suppose that $\text{Out}(\Gamma)$ is infinite. Then $\Gamma$ admits a nontrivial splitting over a finitely generated free abelian group.
\end{theorem}

The results of \cite{CG} now imply the following

\begin{lemma} [cf. Corollary 5.7, \cite{CWIF}] \label{LimitCSA}
Suppose that $\Gamma$ is a torsion-free group hyperbolic relative to a collection of free abelian subgroups, and suppose that $L$ is a $\Gamma$-limit group.  Then
\begin{enumerate}
\item Any finitely generated subgroup of $L$ is a $\Gamma$-limit group;
\item $L$ is torsion-free;
\item $L$ is commutative-transitive and CSA; and
\item Every solvable subgroup of $L$ is abelian.
\end{enumerate}
\end{lemma}

\section{The shortening argument} \label{ShortenSection}

\begin{definition} [Dehn twists]
Let $G$ be a finitely generated group.  A {\em Dehn twist} is an automorphism of one of the following two types:
\begin{enumerate}
\item Suppose that $G = A \ast_C B$ and that $c$ is contained in the centre of $C$.  Then define $\phi \in \text{Aut}(G)$ by $\phi(a) = a$ for $a \in A$ and $\phi(b) = cbc^{-1}$ for $b \in B$;
\item Suppose that $G = A \ast_C$, that $c$ is in the centre of $C$, and that $t$ is the stable letter of this HNN extension.  Then define $\phi \in \text{Aut}(G)$ by $\phi(a) = a$ for $a \in A$ and $\phi(t) = tc$.
\end{enumerate}
\end{definition}

\begin{definition} [Generalised Dehn twists]
Suppose $G$ has a graph of groups decomposition with abelian edge groups, and $A$ is an abelian vertex group in this decomposition.  Let $A_1 \le A$ be the subgroup generated by all edge groups connecting $A$ to other vertex groups in the decomposition.  Any automorphism of $A$ that fixes $A_1$ elementwise can be naturally extended to an automorphism of the ambient group $G$.  Such an automorphism is called a {\em generalised Dehn twist} of $G$.
\end{definition}

\begin{definition} \label{Mod}
Let $G$ be a finitely generated group.  We define ${\rm Mod}(G)$ to be the subgroup of ${\rm Aut}(G)$ generated by:
\begin{enumerate}
\item Inner automorphisms;
\item Dehn twists arising from splittings of $G$ with abelian edge groups; and
\item Generalised Dehn twists arising from graph of groups decompositions of $G$ with abelian edge groups.
\end{enumerate}
\end{definition}
Similar definitions are made in \cite[\S 5]{Sela1} and \cite[\S 1]{BFSela}.

Suppose that $\Gamma$ is a torsion-free group which is hyperbolic relative to abelian subgroups,
acting by isometries on the space $X$ constructed in Section \ref{ConstructX}, with basepoint $x \in X$. Suppose also that $G$ is a finitely generated group, with finite generating set $\A$.  Let $h : G \to \Gamma$ be a homomorphism.  Recall that in Section \ref{Asymp} we defined the {\em length} of $h$ by
\[	\| h \| := \max_{g \in \A} \big\{ d_X(x, h(g) . x) \big\} . 	\]

\begin{definition} [cf. Definition 4.2, \cite{BFSela}] \label{ShortHomo}
We define an equivalence relation on the set of homomorphisms $h : G \to \Gamma$ by setting $h_1 \sim h_2$ if there is $\alpha \in \text{Mod}(G)$ and $\gamma \in \Gamma$ so that $h_1 = \tau_\gamma \circ h_2 \circ \alpha$, where $\tau_\gamma$ is the inner automorphism of $\Gamma$ induced by $\gamma$.

A homomorphism $h : G \to \Gamma$ is {\em short} if for any $h'$ such that $h \sim h'$ we have $\| h \| \le \| h' \|$.
\end{definition}

The following is one of the main technical results of this paper.

\begin{theorem} [Shortening Argument] \label{ShorteningArgument}
Suppose that $\Gamma$ is a non-abelian, freely indecomposable, torsion-free group which is hyperbolic relative to abelian subgroups, and suppose that the sequence of automorphisms $\{ h_n : \Gamma \to \Gamma \}$ converges to a faithful action $\eta : \Gamma \to \text{Isom}(\mathcal C_\infty)$ as above.  Then for all but finitely many $n$ the homomorphism $h_n$ is not short.
\end{theorem}

In order to `shorten' arbitrary homomorphisms, rather than just automorphisms, we need to introduce to new `bending' moves. This is undertaken in \cite{MR-RH} (using ideas inspired
by Alibegovi\'c \cite{Alibegovic2}.

We now show how Theorem \ref{ShorteningArgument} implies Theorem \ref{ModinAut}.

\begin{proof} [Proof (of Theorem \ref{ModinAut}, assuming Theorem \ref{ShorteningArgument}).]
If $\Gamma$ is abelian then the theorem is clear, since in this case ${\rm Mod}(\Gamma) = {\rm Aut}(\Gamma)$.  Thus we assume that $\Gamma$ is non-abelian.  

For each coset $C_i = \rho_i \text{Mod}(\Gamma)$ of $\text{Mod}(\Gamma)$ in $\text{Aut}(\Gamma)$ choose a representative $\hat{\rho}_i$ which is shortest amongst all representatives of $C_i$.  That is to say, each of the automorphisms $\hat{\rho}_i$ is short. 

However, by Theorem \ref{ShorteningArgument} we cannot have an infinite sequence $\{ \hat\rho_n : \Gamma \to \Gamma \}$ of non-equivalent short automorphisms, since then some subsequence will converge to a faithful action of $\Gamma$ on a space $\mathcal C_\infty$.  Hence $\text{Mod}(\Gamma)$ has finite index in $\text{Aut}(\Gamma)$ as required.
\end{proof}

The remainder of this paper is devoted to proving Theorem \ref{ShorteningArgument}.
Before launching into the proof of Theorem \ref{ShorteningArgument}, we need to recall some
of the theory of groups acting on $\R$-trees.

\section{Isometric actions on $\R$-trees} \label{ActionsonRtrees}

In this section we recall a result of Sela from \cite{SelaAcyl}.  Given a finitely generated group $G$ and an $\R$-tree $T$ with an isometric $G$-action, Theorem \ref{RtreeStructure} below gives a decomposition of $T$ which induces a graph of groups decomposition of $G$.  In the case that $G$ is finitely presented, this result follows immediately from Rips Theory; see Bestvina and Feighn, \cite{BF}.

There are two sets of terminology in English for the components of the above-mentioned decomposition\footnote{There is also the work in French by Gaboriau, Levitt and Paulin, with its attendant French terminology;  see \cite{PaulinBour}, for example.}.  Since we are quoting Sela's result, we use his (Rips') terminology.  However, we assume that all axial components are isometric to a real line.  Using Rips and Sela's definition of axial (see \cite[\S 10]{RipsSelaGAFA}), one other case could arise in the arguments that follow (where our group splits as $A \ast_{[a,b]}\langle a,b \rangle$).  Just as noted in \cite[\S4, p.346]{RipsSelaGAFA}, we can treat this case as an IET component.  Thus, without further mention, we consider all axial components to be isometric to a real line.

The following theorem of Sela is used to decompose our limiting $\R$-trees.

\begin{theorem} [Theorem 3.1, \cite{SelaAcyl}; see also Theorem 1.5, \cite{Sela1}] \label{RtreeStructure}
Let $G$ be a freely indecomposable finitely generated group which admits a stable isometric action on a real tree $Y$.  Assume that the stabiliser in $G$ of each tripod in $T$ is trivial.
\begin{enumerate}
\item[1)] There exist canonical orbits of subtrees of $T$, denoted $T_1, \ldots , T_k$, with the following properties:
\begin{enumerate}
\item[(i)] for each $g \in G$ and each $i, j \in \{ 1, \ldots , k \}$ with $i \neq j$, the subtree $g.T_i$ intersects $T_j$ in at most a single point;
\item[(ii)] for each $g \in G$ and each $i \in \{ 1, \ldots , k \}$, the subtree $g.T_i$ is either equal to $T_i$ or intersects $T_i$ in at most a single point;
\item[(iii)] The action of $\text{Stab}_G(T_i)$ on $T_i$ is of axial or IET type;
\end{enumerate}
\item[2)] $G$ is the fundamental group of a graph of groups with:
\begin{enumerate}
\item[(i)] Vertices corresponding to orbits of branching points with non-trivial stabliser in $T$;
\item[(ii)] Vertices corresponding to the orbits of the canonical subtrees $T_1, \ldots , T_k$ which are of axial or IET type.  The groups associated with these vertices are conjugates of the stabilisers of these subtrees.  To a stabiliser of an IET component there exists an associated $2$-orbifold, $\mathcal O$.  Any element of $\pi_1(\mathcal O)$ which corresponds to a boundary component or branching point in $\mathcal O$ stabilises a point in $T$.  For each stabiliser of an IET subtree we add edges that connect the vertex stabilised by it and the vertices stabilised by its boundary components and branching points;
\item[(iii)] Edges corresponding to orbits of edges between branching points with non-trivial stabiliser in the discrete part of $T$ (see Terminology \ref{DiscretePart} below) with edge groups which are conjugates of the stabilisers of these edges;
\item[(iv)] Edges corresponding to orbits of points of intersection between the orbits of $T_1, \ldots , T_k$.
\end{enumerate}
\end{enumerate}
\end{theorem}

\begin{terminology} \label{DiscretePart}
Let $G$ and $T$ be as in Theorem \ref{RtreeStructure} above.  The {\em discrete part of $T$} is the union of the metric closures of the connected components of $T \setminus \left( \bigcup\limits_{i=1}^k G T_i \right)$.
\end{terminology}

\begin{remark} \label{NoThin}
In the theory of stable isometric actions on $\R$-trees, there is one further type of component arising in the decomposition of $T$.  This is called `thin' in \cite{BF} and was discovered and investigated by Levitt (see \cite{Levitt}).  However, in case $G$ is freely indecomposable and the stabiliser of any non-degenerate tripod is trivial (both of these conditions hold in all cases in this paper), thin components do not arise.
\end{remark}

\section{The shortening argument -- Outline} \label{IET}

In this section we outline the proof of Theorem \ref{ShorteningArgument} (the complete proof is contained in this and the subsequent two sections):

\medskip

{\bf Theorem \ref{ShorteningArgument}} (Shortening Argument).
{\em Suppose that $\Gamma$ is a nonabelian, freely indecomposable, torsion-free relatively hyperbolic group with abelian parabolics, and suppose that the sequence of automorphisms $\{ h_n : \Gamma \to \Gamma \}$ converges to an action $\eta : \Gamma \to \text{Isom}(\mathcal C_\infty)$ as above.  Then for all but finitely many $n$ the homomorphism $h_n$ is not short.
}

\medskip

\begin{remark}
Although we call the above theorem the `Shortening Argument', at least for hyperbolic groups the shortening argument is really a collection of ideas applicable in myriad situations.  The above theorem is enough to prove Theorem \ref{ModinAut}.
In order to build Makanin-Razborov diagrams in \cite{MR-RH}, we need to adapt the shortening
argument by adding `bending' moves (see also \cite{Alibegovic2}).
\end{remark}

Let $\{ h_n : \Gamma \to \Gamma \}$ be a sequence of pairwise non-conjugate automorphisms.  Since $\Gamma$ is non-abelian, the action of $\Gamma$ on the limiting space $\mathcal C_\infty$ is faithful, and the action of $\Gamma$ on the associated $\R$-tree $T$ is also faithful.  We prove that for all but finitely many $n$, the homomorphism $h_n$ is not short.

Since the action of $\Gamma$ on $T$ is faithful, $\Gamma$ is itself a strict $\Gamma$-limit group, and by Theorem \ref{LinfProps}.\eqref{tripod} the stabiliser in $\Gamma$ of any tripod in $T$ is trivial.

The approach to proving Theorem \ref{ShorteningArgument} is as follows:  we consider the finite generating set $\A_1$ of $\Gamma$, and the basepoint $y$ of $T$.  We consider the paths $[y, u.y]$ where $u \in \A_1$.  These paths can travel through various types of subtrees of $T$; the IET subtrees, the axial subtrees, and the discrete part of $T$.\footnote{Recall by Remark \ref{NoThin} that there are no thin components in the decomposition of $T$.}  Depending on the types of subtrees which have positive length intersection with $[y,u.y]$, we need various types of arguments which allow us to shorten the homomorphisms which `approximate' the action of $\Gamma$ on $\mathcal C_\infty$.

Mostly, we follow the shortening argument as developed in \cite{RipsSelaGAFA}.  There are two main obstacles to implementing this strategy in the context of torsion-free relatively hyperbolic groups with abelian parabolics.  Note that the automorphisms $h_n : \Gamma \to \Gamma$ actually approximate the action of $G$ on $\mathcal C_\infty$, from which the action of $\Gamma$ on $T$ is extracted.  The two main impediments are: (i) those lines $p_E \in \P$ which correspond to flats $E \in \mathcal C_\infty$; and (ii) that triangles in the approximating spaces are only relatively thin, not actually thin.

\subsection{IET components}

The following theorem of Rips and Sela deals with IET components.

\begin{theorem} \cite[Theorem 5.1, pp. 346-347]{RipsSelaGAFA} \label{IETTheorem}
Let $G$ be a finitely presented freely indecomposable group and assume that $G \times T \to T$ is a small stable action of $G$ on a real tree $T$ with trivial stabilisers of tripods.  Let $U$ be a finite subset of $G$ and let $y \in T$.  Then there exists $\phi_I \in {\rm Mod}(G)$ such that for any $u \in U$, if $[y,u(y)]$ has an intersection of positive length with some {\rm IET}-component of $T$ then:
\[	d_Y(y,\phi_I(u) . y) <d_T(y,u(y)),	\]
and otherwise $\phi_I(u) = u$.
\end{theorem}

It is worth noting that in \cite{RipsSelaGAFA} a more restrictive class of automorphisms is used to shorten the homomorphisms.  Since it {\em is} a more restrictive class, their results still hold using our definition of $\text{Mod}(G)$.

\begin{proposition} \label{ShortenIET}
Suppose that $Y$ is an {\rm IET} subtree of $T$ and that $p_E \in \P$ is a line in $T$.  Then the intersection $Y \cap p_E$ contains at most a point.
\end{proposition}
\begin{proof}
Since $Y$ is an IET subtree, if $\sigma$ is a nondegenerate arc in $Y$ and $\epsilon > 0$ then there exists $\gamma \in \text{Stab}(Y)$ so that $\gamma . \sigma \cap \sigma$ has positive length and there is some $x \in \sigma$ such that $d_T(x, \gamma .x) < \epsilon$.  

Suppose that $Y \cap p_E$ contains more than a point.  By the above remark there exists $\gamma \in \text{Stab}(Y)$ for which $\gamma . p_E \cap p_E$ contains more than a point.  Hence $\gamma . p_E = p_E$, and $p_E \subset Y$.  This, combined with the above fact about IET components, implies that the action of $\text{Stab}(p_E)$ on $p_E$ is indiscrete.  However, this implies that it contains a noncyclic free abelian group, which cannot be a subgroup of $\text{Stab}(Y)$ when $Y$ is an IET subtree.  This contradiction proves the proposition.
\end{proof}

\begin{corollary} \label{IETinCinf}
Let $T$ be an $\R$-tree arising from some $\mathcal C_\infty$ as above.  Suppose $Y$ is an {\rm IET} subtree of $T$ and $\sigma \subset Y$ is a nondegenerate segment.  Then there is a segment $\hat{\sigma} \subset \mathcal C_\infty$, of the same length as $\sigma$, which corresponds to $\sigma$ under the projection from $\mathcal C_\infty$ to $T$.
\end{corollary}

\subsection{Non-IET subtrees, technical results, and the proof of Theorem \ref{ShorteningArgument}}

\ \par

An entirely analogous argument to that of Proposition \ref{ShortenIET} proves

\begin{proposition}
Suppose that a line $l \subset T$ is an axial subtree and the line $p_E \subset T$ is associated to a flat $E \subset \mathcal C_\infty$.  If $l \cap p_E$ contains more than a point then $l = p_E$.
\end{proposition}

\begin{corollary} \label{NotPAxialinCinf}
Let $T$ be an $\R$-tree arising from some $\mathcal C_\infty$ as above.  Suppose $l$ is an axial component of $T$ so that $l \not\in \P$ and $\sigma \subset l$ is a non-degenerate segment.  Then there is a segment $\hat\sigma \subset \mathcal C_\infty$, of the same length as $\sigma$, which corresponds to $\sigma$ under the projection from $\mathcal C_\infty$ to $T$.
\end{corollary}

\begin{lemma} \label{DiscreteEdgeinpE}
If an edge $e$ in the discrete part of $T$ has an intersection of positive length with some line $p_E$ then $e \subset p_E$.
\end{lemma}
\begin{proof}
Suppose that $e$ contains a nontrivial segment from $p_E$ but that $e \not\subset p_E$.  Let $C$ be the edge stabiliser of $e$.  Since $\Gamma$ is freely indecomposable, $C$ is non-trivial, and since $\Gamma$ is torsion-free, $C$ is infinite.  Let $\gamma \in C$.  Then $\gamma$ leaves more than one point of $p_E$ invariant, so leaves all of $p_E$ invariant.  Thus $\gamma$ leaves $E \subset \mathcal C_\infty$ invariant.  Also, since $e \not\subset p_E$, $\gamma$ leaves some point $v \in \mathcal C_\infty \smallsetminus E$ invariant.

By Proposition \ref{FlatInv}, if $\{ E_i \}$ is a sequence of flats ($E_i \subset X_i$) which converges to $E$, then for all but finitely many $n$ the element $h_n(\gamma)$ leaves $E_n$ invariant.  By choosing an $n$ large enough, $h_n(\gamma).E_n = E_n$, and furthermore if $\{ v_i \}$ represents $v$, then $h_n(\gamma)$ moves $v_n$ a distance which is much smaller than the distance from $v_n$ to $E_n$.  In particular, we can ensure that the geodesic $[v_n, h_n(\gamma).v_n]$ does not intersect the $4\delta$-neighbourhood of $E_n$.  Then by Proposition \ref{ProjectProp}, if $\pi : X_n \to E_n$ is the projection map then $d_{X}(\pi (v_n), \pi(h_n(\gamma).v_n)) \le N_3(\phi(3\delta) + N_1(\phi(\delta)))$.  However, since $h_n(\gamma) . E_n = E_n$, we know that $\pi(h_n(\gamma) . v_n) = h_n(\gamma). \pi(v_n)$. Thus $h_n(\gamma)$ moves $\pi(v_n)$ a distance at most $N_3(\phi(3\delta) + N_1(\phi(\delta)))$.

Repeating this argument with a large enough subset of $C$ (namely a subset larger than the maximal size of an intersection of any orbit $\Gamma . u$ with a ball of radius $N_3(\phi(3\delta) + N_1(\phi(\delta)))$), we obtain a (finite) bound on the size of $C$.  However, $C$ is infinite, as noted above.  This contradiction finishes the proof.
\end{proof}

The following Theorems \ref{AxialTheorem}, \ref{pETheorem} and \ref{DiscretewithpE} are the technical results needed to prove Theorem \ref{ShorteningArgument}.

\medskip

{\bf Theorem \ref{AxialTheorem}.}
{\em Let $G$ be a finitely generated freely indecomposable group and assume that $G \times T \to T$ is a small stable action of $G$ on an $\R$-tree $T$ with trivial stabilisers of tripods.  Let $U$ be a finite subset of $G$ and let $y \in T$.  Then there exists $\phi_A \in {\rm Mod}(G)$ so that for any $u \in U$, if $[y, u.y]$ has an intersection of positive length with some axial component of $T$ then:
\[	d_T(y, \phi_A(u) . y) < d_T(y, u.y) ,	\]
and otherwise $\phi_A(u) = u$.
}

\medskip

As far as I am aware, Theorem \ref{AxialTheorem} has not appeared in print.  However, its statement and proof are very similar to those of Theorem \ref{IETTheorem}, and it is certainly known at least to Sela (see \cite[\S 5]{Sela1}) and to Bestvina and Feighn (see \cite[Exercise 11]{BFSela}).

\begin{remark} \label{fgVSfp}
Theorem \ref{IETTheorem} is stated for finitely presented groups.  The only time in the proof when it is required that $G$ be finitely presented rather than just finitely generated is when a result of Morgan from \cite{Morgan} is quoted.

Specifically, Rips and Sela show that when $G$ is freely indecomposable and finitely presented and acts on an $\R$-tree $T$ with trivial tripod stabilisers then, for each $g \in G$ and any $y \in T$, the path $[y,\gamma . y]$ cuts only finitely many components of axial or IET type in $Y$ (see \cite[pp. 350--351]{RipsSelaGAFA}).  

However, this is also true when $G$ is only assumed to be finitely generated, rather than finitely presented (but all other assumptions apply).  This follows from the arguments in \cite[\S 3]{SelaAcyl}.  The action of $G$ on $T$ can be approximated by actions of finitely presented groups $G_i$ on $\R$-trees $Y_i$.  For large enough $k$, the IET and axial components of $Y_k$ map isometrically onto the IET and axial components of $T$ (see \cite[\S 3]{SelaAcyl} for details).

Therefore, Theorem \ref{IETTheorem} stills holds when $G$ is assumed to be finitely generated, but not necessarily finitely presented.  Similarly, Theorem \ref{AxialTheorem} above, whose proof mimics the proof of Theorem \ref{IETTheorem}, holds finitely generated groups $G$.  However, in this paper we can assume $G$ is finitely presented.
\end{remark}

We now state the further technical results which are required for the proof of Theorem \ref{ShorteningArgument}.  These technical results are proved in the subsequent two sections.

\medskip

{\bf Theorem \ref{pETheorem}.}
{\em Let $\Gamma$ be a freely indecomposable, torsion-free, non-abelian relatively hyperbolic group with abelian parabolics.  Suppose that $h_n : \Gamma \to \Gamma$ is a sequence of automorphisms converging to a faithful action of $\Gamma$ on a limiting space $\mathcal C_\infty$ and let $T$ be the $\R$-tree associated to $\mathcal C_\infty$.  Let $U$ be a finite subset of $\Gamma$.  Let $y \in T$, let $\hat{y} \in \mathcal C_\infty$ project to $y \in T$ and let $\{ \hat{y}_m \}$ be a sequence representing $\hat{y}$.  Let $p_E$ be an axial component of $T$.  There exists $m_0$ so that: for all $m \ge m_0$ there is $\phi_{p_E,m} \in {\rm Mod}(\Gamma)$ so that for any $u \in U$, if $[y,u.y]$ has an intersection of positive length with a line in the $\Gamma$-orbit of $p_E$ then
\[	d_{X_m}(\hat{y}_m, (h_m \circ \phi_{p_E,m})(u)) . \hat{y}_m) < d_{X_m}(\hat{y}_m,h_m(u). \hat{y}_m) ,	\]
and otherwise $\phi_{p_E,m}(u) = u$.
}

\medskip

{\bf Theorem \ref{DiscretewithpE}.}
{\em Let $\Gamma$ be a freely indecomposable torsion-free relatively hyperbolic group with abelian parabolics.  Suppose that $h_n : \Gamma \to \Gamma$ is a sequence of automorphisms converging to a faithful action $\Gamma$ on a limiting space $\mathcal C_\infty$ with associated $\R$-tree $T$.  Suppose further that $U$ is a finite subset of $\Gamma$.  Let $y \in T$, let $\hat{y} \in \mathcal C_\infty$ project to $y \in T$ and let $\{ \hat{y}_m \}$ be a sequence representing $\hat{y}$.  There exists $m_0$ so that: for all $m \ge m_0$ there is $\phi_{D,m} \in {\rm Mod}(\Gamma)$ so that for any $u \in U$ which does not fix $y$ and with $[y,u.y]$ supported only in the discrete parts of $T$ we have
\[	d_{X_m}(\hat{y}_m, (h_m \circ \phi_{D,m})(u) . \hat{y}_m) < d_{X_m}(\hat{y}_m, u .\hat{y}_m) .	\]
}

\medskip

Armed with Theorem \ref{ShortenIET}, and assuming Theorems \ref{AxialTheorem}, \ref{pETheorem} and \ref{DiscretewithpE}, we now prove Theorem \ref{ShorteningArgument}.

\begin{proof}[Proof (Theorem \ref{ShorteningArgument}).]
We have already noted that the action of $\Gamma$ on $T$ is faithful, that $T$ is not isometric to a real line, and that the stabiliser in $\Gamma$ of any tripod in $T$ is trivial.  Also, $\SK (h_n) = \{ 1 \}$.  

We suppose (by passing to a subtree if necessary) that the tree $T$ is minimal. As noted in Remark \ref{NoThin} above, $T$ contains no thin components.

Let $U = \A_1$ be the fixed generating set of $\Gamma$ used to define $\| f \|$ for a homomorphism $f : \Gamma \to \Gamma$, let $y$ be the image in $T$ of the basepoint $x_\omega \in \mathcal C_\infty$ and let $\hat{y}_m = x$ for each $m$.

Let $\phi_I$ be the automorphism of $\Gamma$ given by Theorem \ref{IETTheorem} and $\phi_A$ the automorphism from Theorem \ref{AxialTheorem}.  

Suppose that $u \in U$ is such that $[y,u.y]$ has an intersection of positive length with an IET component of $T$.  Then Theorem \ref{IETTheorem} and Corollary \ref{IETinCinf} guarantee that for all but finitely many $n$ we have $\| h_n \circ \phi_I \| < \| h_n \|$ so $h_n$ is not short.  Similarly, if $[y,u.y]$ has an intersection of positive length with an axial component which is not contained in any $p_E \in \mathcal P$ then Theorem \ref{AxialTheorem} and Corollary \ref{NotPAxialinCinf} imply that for all but finitely many $n$ we have $\| h_n \circ \phi_A \| < \| h_n \|$ so also in this case $h_n$ is not short.

Suppose then that $[y,u.y]$ has an intersection of positive length with a line in the $\Gamma$-orbit of some $p_E$, and suppose that $p_E$ is an axial component of $T$.  Then by Theorem \ref{pETheorem} for all but finitely many $n$ there exists an automorphism $\phi_{p_E,n} \in \text{Mod}(\Gamma)$ so that $\| h_n \circ \phi_{p_E,n} \| < \| h_n \|$, so $h_n$ is not short.

Finally, suppose that all of the segments $[y, u.y]$ are entirely contained in the discrete part of $T$.  Then by Theorem \ref{DiscretewithpE} for all but finitely many $n$ there exists $\phi_{D,n} \in {\rm Mod}(\Gamma)$ so that $\| h_n \circ \phi_{D,n} \| < \| h_n \|$, and once again $h_n$ is not short.

This completes the proof of the theorem.
\end{proof}

Having proved Theorem \ref{ShorteningArgument} we now prove Theorem \ref{ModinAut}.  Given Theorem \ref{ShorteningArgument}, the proof is identical to that of \cite[Corollary 4.4]{RipsSelaGAFA}.

\medskip

{\bf Theorem \ref{ModinAut}}
{\em Suppose that $\Gamma$ is a freely indecomposable torsion-free relatively hyperbolic group with abelian parabolics.  Then ${\rm Mod}(\Gamma)$ has finite index in ${\rm Aut}(\Gamma)$.}

\smallskip

\begin{proof}
If $\Gamma$ is abelian then the theorem is clear, since in this case ${\rm Mod}(\Gamma) = {\rm Aut}(\Gamma)$.  Thus we assume that $\Gamma$ is non-abelian.  

For each coset $C_i = \rho_i \text{Mod}(\Gamma)$ of $\text{Mod}(\Gamma)$ in $\text{Aut}(\Gamma)$ choose a representative $\hat{\rho}_i$ which is shortest amongst all representatives of $C_i$.  That is to say, each of the automorphisms $\hat{\rho}_i$ is short. 

However, by Theorem \ref{ShorteningArgument} we cannot have an infinite sequence $\{ \hat\rho_n : \Gamma \to \Gamma \}$ of non-equivalent short automorphisms, since then some subsequence will converge to a faithful action of $\Gamma$ on a space $\mathcal C_\infty$.  Hence $\text{Mod}(\Gamma)$ has finite index in $\text{Aut}(\Gamma)$ as required.
\end{proof}

\section{Axial components} \label{AxialSection}

The purpose of this section is to prove the following two theorems.

\begin{theorem} \label{AxialTheorem}
Let $G$ be a finitely generated freely indecomposable group and assume that $G \times T \to T$ is a small stable action of $G$ on an $\R$-tree $T$ with trivial stabilisers of tripods.  Let $U$ be a finite subset of $G$ and let $y \in T$.  Then there exists $\phi_A \in {\rm Mod}(G)$ so that for any $u \in U$, if $[y, u.y]$ has an intersection of positive length with some axial component of $T$ then:
\[	d_T(y, \phi_A(u) . y) < d_T(y, u.y) ,	\]
and otherwise $\phi_A(u) = u$.
\end{theorem}

\begin{theorem} \label{pETheorem}
Let $\Gamma$ be a freely indecomposable, torsion-free, non-abelian relatively hyperbolic group with abelian parabolics.  Suppose that $h_n : \Gamma \to \Gamma$ is a sequence of automorphisms converging to a faithful action of $\Gamma$ on a limiting space $\mathcal C_\infty$ and let $T$ be the $\R$-tree associated to $\mathcal C_\infty$.  Let $U$ be a finite subset of $\Gamma$.  Let $y \in T$, let $\hat{y} \in \mathcal C_\infty$ project to $y \in T$ and let $\{ \hat{y}_m \}$ be a sequence representing $\hat{y}$.  Let $p_E$ be an axial component of $T$.  There exists $m_0$ so that: for all $m \ge m_0$ there is $\phi_{p_E,m} \in {\rm Mod}(\Gamma)$ so that for any $u \in U$, if $[y,u.y]$ has an intersection of positive length with a line in the $\Gamma$-orbit of $p_E$ then
\[	d_{X_m}(\hat{y}_m, (h_m \circ \phi_{p_E,m})(u)) . \hat{y}_m) < d_{X_m}(\hat{y}_m,h_m(u). \hat{y}_m) ,	\]
and otherwise $\phi_{p_E,m}(u) = u$.
\end{theorem}

To prove Theorem \ref{AxialTheorem} we follow the proof of \cite[Theorem 5.1]{RipsSelaGAFA} (which is Theorem \ref{IETTheorem} in this paper).  First, we need the following result, the (elementary) proof of which we include because of its similarity to Proposition \ref{ShortenpEAxial} below.

\begin{proposition} \label{ShortenAxial}
Suppose that $\rho: P \times \R \to \R$ is an orientation-preserving, indiscrete isometric action of $P \cong \Z^n$ on the real line $\R$.  For any finite subset $W$ of $P$and any $\epsilon > 0$ there exists an automorphism $\sigma : P \to P$ such that:
\begin{enumerate}
\item[1)] For every $w \in W$ and every $r \in \R$
\[	d_{\R}(r, \sigma (w) . r ) < \epsilon ;	\]
\item[2)] For any $k \in \text{ker}(\rho)$ we have $\sigma(k) = k$.
\end{enumerate}
\end{proposition}
\begin{proof}
There is a direct product decomposition $P = A \oplus B$ where $A$ is the kernel of the action of $P$ on $\R$, and $B$ is a finitely generated free abelian group which has a free, indiscrete and orientation preserving action on $\R$.  The automorphism $\sigma$ we define fixes $A$ elementwise, so we can assume that all elements of $W$ lie in $B$ (since elements of $A$ fix $\R$ pointwise).  Thus, we need only prove the lemma in case the action is faithful.

Since the action of $B$ on $\R$ is indiscrete and free, the translation lengths of elements of a basis of $B$ are independent over $\Z$.  In particular, there is a longest translation length amongst the translation lengths of a basis of $B$.  Suppose that $b_1 \in B$ is the element of the basis with largest translation length, and that $b_2$ has the second largest.  Denote these translation lengths by $| b_1 |$ and $| b_2 |$, respectively.  Since $| b_1 |$ and $| b_2 |$ are independent over $\Z$, there is $n \in \Z$ so that $0 < |b_1 + nb_2| < | b_2 |$.  Replace $b_1$ by $b_1 + nb_2$.  This is an automorphism of $P$, fixing $A$ elementwise.  

Proceeding in this manner, we can make each of the elements of a basis as small as we wish, and so given $W$ and $\epsilon > 0$, we can make each of the elements of $W$ (considered as a word in the basis of $B$) have translation length less than $\epsilon$, as required.
\end{proof}

\begin{proof}[Proof (Theorem \ref{AxialTheorem}).]

By Claim \ref{fgVSfp}, each of the segments $[y,u.y]$ for $u \in U$ cuts only finitely many components of $T$ of axial or IET type.  Let $\epsilon$ be the minimum length of a (non-degenerate) interval of intersection between $[y,u.y]$ and an axial component of $T$, for all $u \in U$.

The action of $G$ on $T$ induces a graph of groups decomposition $\Lambda$ of $G$ as in Theorem \ref{RtreeStructure}.  Let $T_i$ be an axial component of $T$.  There is a vertex group of $\Lambda$ corresponding to the $G$-orbit of $T_i$, with vertex group a conjugate of $\text{Stab}(T_i)$.  By Theorem \ref{LinfProps} and Lemma \ref{LimitCSA}, $\text{Stab}(T_i)$ is a free abelian group.  The vertex groups adjacent to $\text{Stab}(T_i)$ (those that are separated by a single edge) stabilise a point in the orbit of a branching point in $T_i$ with nontrivial stabiliser.  Recall that $G$ is freely indecomposable, so all edge groups are nontrivial.

Let $q_1$ be the point on $T_i$ closest to $y$ (if $y \in T_i$ then $q_1 = y$).  Choose points $q_2, \ldots , q_m \in T_i$ in the orbits of the branching points corresponding to the adjacent vertex groups such that $d_T(q_i,q_j) < \frac{\epsilon}{20}$.  We can do this since the action of $\text{Stab}(T_i)$ on $T_i$ has all orbits dense, since $T_i$ is an axial component.

The proof of Theorem \ref{AxialTheorem} now proceeds exactly as the proof of \cite[Theorem 5.1]{RipsSelaGAFA} (start with Case 1 on p.351).
\end{proof}

\begin{proof} [Proof (Theorem \ref{pETheorem}).]
Since $p_E \in \P$ is an axial component of $T$, there is a vertex group corresponding to the conjugacy class of $\text{Stab}(p_E)$ in the graph of groups decomposition which the (faithful) action of $\Gamma$ on $T$ induces (see Theorem \ref{RtreeStructure}).

Now, the stabiliser in $\Gamma$ of $p_E$ is exactly the stabiliser in $\Gamma$ of $E$, when $\Gamma$ acts (also faithfully) on $\mathcal C_\infty$.  By \cite[Corollary 3.17]{CWIF}, there is a sequence of flats $E_i$ in the approximating spaces $X_i$ so that $E_i \to E$ in the $\Gamma$-equivariant Gromov topology. By Proposition \ref{FlatInv}, if $\gamma \in \text{Stab}_\Gamma(E)$ then for all but finitely $i$ we have $h_i(\gamma) \in \text{Stab}(E_i)$. For such an $i$, the element $h_i(\gamma)$ is contained in a unique noncyclic maximal abelian subgroup $A_i$  of $\Gamma$.  
However, $h_i$ is an automorphism, so $\gamma$ is contained in a unique noncyclic maximal abelian subgroup $A_\gamma$ of $\Gamma$, and $A_i = h_i(A_\gamma)$.

If $\gamma'$ is another element of $\text{Stab}_\Gamma(E)$, then $[\gamma, \gamma'] = 1$, and it is not difficult to see that $A_\gamma = A_{\gamma'}$.  Also, if $\gamma_0 \in A_\gamma$ then $\gamma_0 \in \text{Stab}_\Gamma(E)$.  Hence $A_\gamma = \text{Stab}_\Gamma(E)$. We denote the subgroup $\text{Stab}_\Gamma(E)$ by $A_E$.

We now prove Theorem \ref{pETheorem} by finding an analogue of Proposition \ref{ShortenAxial} in the flats $E_i$ and then once again following the proof from \cite{RipsSelaGAFA}.

\begin{proposition} \label{ShortenpEAxial}
Let $W$ be a finite subset of $A_E$.  For any $\epsilon > 0$ there exists $i_0$ so that for all $i \ge i_0$, there is an automorphism $\sigma_i : A_E \to A_E$ so that 
\begin{enumerate}
\item For every $w \in W$, and every $r_i \in E_i$,
\[	d_{X_i}(r_i, h_i(\sigma_i(w)) . r_i) < \epsilon ;	\]
\item For any $k \in A_E$ which acts trivially on $E$ we have $\sigma_i(k) = k$.
\end{enumerate}
\end{proposition}
\begin{proof}[Proof (Proposition \ref{ShortenpEAxial}).]
The group $A_E$ admits a decomposition $A_E = A_0 \oplus A_1$, where $A_0$ acts trivially on $E$, and $A_1$ acts freely on $E$.  Choose a basis $\B$ of $A_E$ consisting of a basis for $A_0$ and a basis for $A_1$.  Let $k_W$ be the maximum word length of any element of $W$ with respect to the chosen basis.

Since the $h_i : \Gamma \to \Gamma$ are automorphisms, for sufficiently large $i$ and any $a \in E_i$, the set $h_i(A_E) . a \subset E_i$ forms an $\frac{\epsilon}{20k_W}$-net in $E_i$ (where distance is measured in the metric $\frac{1}{\| h_i \|}$ on $X_i$).  Choose a (possibly larger) $i$ so that also the action of $h_i(\B)$ on $E_i$ approximates the action of $\B$ on $E$ to within $\frac{\epsilon}{20k_W}$ (note that since the action of $A_E$ on $E$ and the action of $h_i(A_E)$ on $E_i$ are both by translations, and translations of Euclidean space move every point the same distance, there are arbitrarily good approximations for the action of any finite subset of $A_E$ on the whole of $E$).

The remainder of the proof proceeds just as the proof of Propostion \ref{ShortenAxial} above, although in the step where we replace $b_1$ by $b_1 + nb_2$, we cannot insist that $b_2$ acts nontrivially on $E$.  However, we of course can insist that $b_1$ acts nontrivially on $E$, since otherwise it moves all points of $E$ a distance at most $\frac{\epsilon}{20}$.  Therefore, such an automorphism is nonetheless a generalised Dehn twist.
\end{proof}
Given Proposition \ref{ShortenpEAxial}, the proof of Theorem \ref{pETheorem} once again follows the proof of \cite[Theorem 5.1, pp. 350--353]{RipsSelaGAFA}, although in this case we have to choose approximations to the action of $\Gamma$ on $\mathcal C_\infty$ (the important point here is that the sets $h_i(A_E) . a$, for any $a \in E_i$, get denser and denser in $E_i$, when considered in the scaled metric $\frac{1}{\| h_i \|}d_X$ of $X_i$).  These small changes are straightforward, but do lead to the different shortening automorphisms $\phi_{p_E,m}$ in the statement of Theorem \ref{pETheorem}.
\end{proof}

\section{The discrete case} \label{Discrete}

In this section we shorten the approximations to paths of the form $[\hat{y},u. \hat{y}]$, where $\hat{y} \in \mathcal C_\infty$ projects to $y \in T$ and $[y,u.y]$ is entirely supported in the discrete part of $T$.  The lengths of the limiting paths $[\hat{y}, u. \hat{y}]$ and $[y,u.y]$ are unchanged.

The purpose of this section is to prove the following

\begin{theorem} \label{DiscretewithpE}
Let $\Gamma$ be a freely indecomposable torsion-free relatively hyperbolic group with abelian parabolics.  Suppose that $h_n : \Gamma \to \Gamma$ is a sequence of automorphisms converging to a faithful action $\Gamma$ on a limiting space $\mathcal C_\infty$ with associated $\R$-tree $T$.  Suppose further that $U$ is a finite subset of $\Gamma$.  Let $y \in T$, let $\hat{y} \in \mathcal C_\infty$ project to $y \in T$ and let $\{ \hat{y}_m \}$ be a sequence representing $\hat{y}$.  There exists $m_0$ so that: for all $m \ge m_0$ there is $\phi_{D,m} \in {\rm Mod}(\Gamma)$ so that for any $u \in U$ which does not fix $y$ and with $[y,u.y]$ supported only in the discrete parts of $T$ we have
\[	d_{X_m}(\hat{y}_m, (h_m \circ \phi_{D,m})(u) . \hat{y}_m) < d_{X_m}(\hat{y}_m, u .\hat{y}_m) .	\]
\end{theorem}

The proof of Theorem \ref{DiscretewithpE} follows \cite[\S 6]{RipsSelaGAFA}.

By Lemma \ref{DiscreteEdgeinpE}, if $e$ is a discrete edge in $T$ then either $e \in p_E$ for some flat $E \subset \mathcal C_\infty$, or $\mathcal C_\infty$ contains a well-defined, canonical, isometric image $\hat{e}$ of $e$, so that $\hat{e}$ projects to $e$.

We have a sequence of automorphisms $\{ h_n : \Gamma \to \Gamma \}$, converging to a faithful action of $\Gamma$ on a limiting space $\mathcal C_\infty$, with associated $\R$-tree $T$. 

There are a number of different cases to consider:

\smallskip

\noindent{\bf \underline{Case 1:}}  $y$ is contained in the interior of an edge $e$

\smallskip

\noindent{\bf Case 1a:}  $e$ is not completely contained in a line of the form $p_E$ and $\bar{e} \in T/\Gamma$ is a splitting edge.  

Note that because $e$ is not contained in any $p_E$, there is a single point $\hat{y} \in \mathcal C_\infty$ which corresponds to $y \in T$.

This case is very similar to the Case 1a on pp. 355--356 of \cite{RipsSelaGAFA}.  In this case we have a decomposition $\Gamma = A \ast_C B$ where $C$ is a finitely generated free abelian group properly contained in both $A$ and $B$.

Given $u \in U$ we can write:
\[	u = a_u^1b_u^1 \cdots a_u^{n_u}b_u^{n_u} ,	\]
where $a_u^i \in A$ and $b_u^i \in B$ (possibly $a_u^1$ and/or $b_u^{n_u}$ are the identity).  Let $\{ z_1, \ldots , z_n \}$ be a generating set for $Z$.

Let $\epsilon$ be the minimum of:
\begin{enumerate}
\item the length of the shortest edge in the discrete part of $T$;
\item the distance between $y$ and the vertices of $e$.
\end{enumerate}

Recall that triangles in $X$ are relatively $\delta$-thin, and the function $\phi$ comes from the definition of isolated flats.  Let $C_0$ be the maximum size of an intersection of an orbit $\Gamma . z$ with a ball of radius $10 \delta + 2\phi(3\delta)$ in $X$ (where distance is measured in $d_X$).

Now take $F$ to be the finite subset of $G$ containing $1$ and
\[	z a_u^i z^{-1},\  \ \ zb_u^iz^{-1} ,	\]
where $z \in C$ has word length at most $10C_0$.

For large enough $m$ we have, for all $\gamma_1 , \gamma_2 \in F$,

\begin{equation} \label{Approx}
| d_{X_m}(h_m(\gamma_1) . \hat{y}_m, h_m(\gamma_2).\hat{y}_m) - d_{\mathcal C_\infty}(\gamma_1 . \hat{y}, \gamma_2 . \hat{y})| < \epsilon_1 ,
\end{equation}
where $\epsilon_1 = \frac{\epsilon}{8000C_0}$.

Let $w_m \in [\hat{y}_m,h_m(a_u^1).\hat{y}_m]$ and $w_m' \in [\hat{y}_m, h_m(b_u^1).\hat{y}_m]$ satisfy 
\begin{equation} \label{Flexible}
\frac{\epsilon}{2} - \frac{\epsilon}{1000} \le d_{X_m}(w_m,\hat{y}_m) = d_{X_m}(w_m',\hat{y}_m) \le \frac{\epsilon}{2} + \frac{\epsilon}{1000}.
\end{equation}

\begin{lemma} \label{zMoves}
For some $z \in C$ of word length at most $10C_0$ we have, for all but finitely many $m$,
\begin{eqnarray*}
d_{X_m}(\hat{y}_m, h_m(z). w_m) & < & d_{X_m}(\hat{y}_m,w_m) - 8\delta_m, \mbox{ and } \\
d_{X_m}(\hat{y}_m, h_m(z).w_m') & < & d_{X_m}(\hat{y}_m,w_m')+ 8\delta_m.
\end{eqnarray*}
\end{lemma}
\begin{proof}
Let $W$ be the set of all elements $z \in C$ of word length at most $10C_0$ in the generators $\{ z_1, \ldots , z_n \}$ and their inverses.

First suppose that for all but finitely many $i$ we have $h_i(W) \subseteq \text{Stab}_\Gamma(E_i)$.  Then since the edge containing $y$ is not completely contained in a single $p_E$, we can assume that each element of $W$ fixes a point outside of $E$.  Now, using Proposition \ref{ProjectProp}, there is a point in $E_i$ which is moved at most $N_3(\phi(3\delta) + N_1(\phi(\delta)))$ by each element of $h_i(W)$.  This gives a bound on the size of $h_i(W)$ which does not depend on $i$ (so long as $i$ is large enough).  However, this contradicts the choice of $W \subseteq \Gamma$.  Therefore it is not the case that $h_i(W) \subseteq \text{Stab}_\Gamma(E_i)$ for all but finitely many $i$.

By the argument in the paragraphs after the proof of \cite[Lemma 4.5]{CWIF}, for all but finitely many $k$, the elements $h_i(z)$ act approximately like translations.  Since $W$ is closed under inverses, and we have chosen $W$ large enough that some element `translates' by at least $10\delta_m$, we can choose some $z \in W$ which satisfies the conclusion of the lemma.
 \end{proof}
 
 In order to finish Case 1a, we follow the proofs of Proposition 6.3 and Theorem 6.4 from \cite{RipsSelaGAFA}.  The only additional thing needed in this case is to force $w_m$ to lie close to each $[\hat{y}_m, h_m(a_u^i).\hat{y}_m]$.  We do this by applying Lemma \ref{4.5} and the arguments in the paragraphs in \cite{CWIF} which follow the proof of \cite[Lemma 4.5]{CWIF}.  It is for this reason that we left some flexibility as to the choice of $w_m$ and $w_m'$ in \ref{Flexible} above.  
 
 Doing this, we can ensure that $w_m$ lies within $2\delta$ of each geodesic segment $[\hat{y}_m, h_m(a_u^i) . \hat{y}_m]$, and similarly for $w_m'$.  We can now follow the proof \cite[Proposition 6.3]{RipsSelaGAFA}.  The proof of \cite[Theorem 6.4]{RipsSelaGAFA} is not included in \cite{RipsSelaGAFA} (or in \cite{SelaAcyl} as claimed in \cite{RipsSelaGAFA}).  However, it is straightforward, so we omit it here also.
 
The automorphism we use to shorten in this case is:
  \begin{eqnarray*}
 \forall a \in A && \phi(a) = zaz^{-1} \\
 \forall b \in B && \phi(b) = z^{-1}bz,
 \end{eqnarray*}
 where $a$ is as in Lemma \ref{zMoves} above.  This completes the proof in Case 1a.  It is worth noting here that we are shortening the actions on $X_i$ which approximate the action on $\mathcal C_\infty$.  However, this does not affect the analogy between the proofs here and those in \cite{RipsSelaGAFA}.

\smallskip

\noindent{\bf Case 1b:} $e$ is not completely contained in a single $p_E$ and $\bar{e} \in T/\Gamma$ is not a splitting edge.  

In this case we have a decomposition $\Gamma = A \ast_C$, where $C$ is a finitely generated free abelian group.  

In the same way as we adapted the proof of Case 1a from \cite{RipsSelaGAFA} above, we may adapt the proof of Case 1b from \cite{RipsSelaGAFA}.  The key point is that we allow a small amount of flexibility in the choice of $w_m$ and $w_m'$.  Doing this, we may ensure that even though the approximating triangles we consider are only relatively thin, rather than actually thin, all of the features we need to apply the proof from \cite{RipsSelaGAFA} still hold, because we can make sure that we are not near the `fat' part of any triangle.  Proceeding with this idea in mind, the proof from \cite{RipsSelaGAFA} can be adapted without difficulty.

\medskip

We now deal with the two cases where $y$ is contained in the interior of the edge $e$ and $e \subset p_E$ for some $p_E \in \P$.  Using Lemma \ref{222} and Proposition \ref{ProjectProp}, the following result is not difficult to prove:

\begin{proposition}
Suppose that $X$ is the space constructed in Section \ref{ConstructX}. There exists a constant $N_4$, depending only on $X$ so that if  $E_1, E_2 \in \mathcal Q$ are maximal flats in $X$ then
there is a set $J_{E_1,E_2}$ so that:
\begin{enumerate}
\item ${\rm Diam}(J_{E_1,E_2}) \le N_4$; and
\item If $x \in E_1$ and $u \in E_2$ then any geodesic between $x$ and $y$ intersects
$J_{E_1,E_2}$ nontrivially.
\end{enumerate}  
\end{proposition}
Recall that $\mathcal Q$ is the family of maximal flats from the definition of $X$, that triangles in $X$ are relatively $\delta$-thin (Theorem \ref{RelThinTriangles}), and that $\phi$ is the function from Lemma \ref{IsolatedFlats}.  Following Convention \ref{IsolConv}, we assume without loss of generality that for all $k \ge 0$ we have $\phi(k) \ge k$ and also that $\phi$ is a nondecreasing function.

Choose compact fundamental domains for the action of $\text{Stab}_\Gamma(E)$ on $E$, for each conjugacy class of maximal flat in $X$, and let $K_F$ be the maximal diameter of these fundamental domains.  Also, let $K_X$ be the diameter of a compact set $D$ for which $\Gamma . D = X$.  For the remainder of Case 1, we replace the constant $\delta$ by
\[	\max \left\{ \delta , 1000K_F , 1000 K_X, 1000(7\delta + 14\phi(4\delta)) \right\}	.	\]

The stabiliser of the edge $e$ is a subgroup of $\text{Stab}_\Gamma(E)$.  Since $p_E$ is not an axial component, the action of $\text{Stab}(E)$ on $E$ is either trivial or factors through a infinite cyclic group.  If $\bar{e} \in T/\Gamma$ is a splitting edge, then necessarily the action of $\text{Stab}(E)$ on $E$ is trivial. 

\noindent{\bf Case 1c:}  $e$ is completely contained in some $p_E$, and $\bar{e} \in T/\Gamma$  is a splitting edge.

Let $A_E = \text{Stab}_\Gamma(E)$.  Then, $A_E = A_0 \oplus A_1$, where $A_0$ acts trivially on $E$ and $A_1$ acts freely on $E$.  Since $p_E$ is a splitting edge, $A_1 = \{ 1 \}$.

We have a decomposition $\Gamma = H_1 \ast_{A_E} H_2$.

The subgroup $H_1$ fixes a point in $p_E$, but does not fix all of $p_E$.  Thus, $H_1$ fixes a point $v_1 \in E$.  Similarly, $H_2$ fixes a point $v_2 \in E$, but does not fix all of $E$.  We choose points $\hat{y}_m \subset E_m$ so that: (i) $\{ \hat{y}_m \}$ represents $\hat{y} \in \mathcal C_\infty$ which projects to $y \in T$; (ii) each $\hat{y}_m$ lies in the orbit $\Gamma . x$; and (iii) subject to the first two conditions, $\hat{y}_m$ lies as close as possible to the line $[v_1^m, v_2^m]$, where $\{ v_i^m \} \to v_i$, $i = 1,2$.

We proceed as in Case 1a. However, this time we cannot find a single automorphism to shorten the $\| h_i \|$, but we use the fact that the sets $h_i(A_E) .a \subset E_i$ are denser and denser (when distance is measured in the metrics $\frac{1}{\| h_i \|}d_X$) to find, for all but finitely many $i$, a Dehn twist $\phi_{e,i}$ which shortens the action on $X_i$.  This proceeds in a similar way to Case 1a above, using the ideas in Proposition \ref{ShortenpEAxial} and the proof of Theorem \ref{pETheorem} above.

\smallskip

\noindent{\bf Case 1d:} $e$ is completely contained in some $p_E$ and $\bar{e} \in T/\Gamma$ is not a splitting edge.

There are two cases here.  As in Case 1b, we have a decomposition $\Gamma = A \ast_C$, where $C$ is a finitely generated free abelian group.  Let $t$ be the stable letter of this HNN extension, and suppose that $C \le \text{Stab}(E)$, a maximal flat in $\mathcal C_\infty$.  The two cases are where $f \in \text{Stab}(E)$, and when $f \not\in \text{Stab}(E)$.

Each of these cases follow the proof of Case 1b above (and therefore Case 1b from \cite{RipsSelaGAFA}) in the same way as Case 1c followed the proof of Case 1a.

\smallskip

\noindent{\bf Case 2:}  $y$ is a vertex of $T$.

In this case, we do not shorten the approximations to a particular edge, but each of the edges adjacent to $y$.  As before, we largely follow \cite[\S 6]{RipsSelaGAFA}.

There are four cases again, when the edge is splitting, and non-splitting, coupled with the cases where the edge is contained in some $p_E$ and when it is not.

These follow the proofs from \cite{RipsSelaGAFA} just as in Case 1 above.  Note that the shortening automorphisms fix elementwise $\text{Stab}_\Gamma(\hat{y})$.

\begin{proof}[Proof (Theorem \ref{DiscretewithpE})]
If $y$ is contained in the interior of an edge, then apply Case 1 above to find a sequence of automorphisms which shorten the $h_n$.

If $y$ is a vertex in $T$, then we shorten the $h_n$ on each of the adjacent edges separately using Case 2 and \cite[\S 6]{RipsSelaGAFA}.
\end{proof}

This finally completes the proof of Theorem \ref{ShorteningArgument}.

\end{document}